\documentclass[11pt,oneside,a4paper]{amsart}
\usepackage{graphicx}
\usepackage{subfigure}
\usepackage{amsmath,amssymb,amsthm,bbm,marvosym}
\usepackage{enumitem}
\usepackage[]{units}
\setenumerate{leftmargin=*}

\usepackage{fullpage}

\let\subset\subseteq 
\let\eps\varepsilon
\let\rho\varrho

\def\NN{\mathbb{N}}

\def\essinf{\mathrm{essinf\,}}

\newcommand{\JUSTIFY}[1]{\fbox{\tiny{#1}}\quad}

\def\differential{\mathsf{d}} 

\def\Pr{\mathbf{P}} 
\def\Exp{\mathbf{E}} 
\def\Bin{\mathsf{Bin}} 
\def\Poi{\mathsf{Poi}} 
\def\RG{\mathbb{G}} 
\def\RK{\mathbb{K}} 
\def\cI{\mathcal{I}}
\def\cJ{\mathcal{J}}
\def\cE{\mathcal{E}}
\def\cD{\mathcal{D}}
\def\bh{\mathbf{h}}

\newtheorem{theorem}{Theorem}
\newtheorem{lemma}[theorem] {Lemma}   
   
\newtheorem{conjecture}[theorem] {Conjecture}  
\newtheorem{proposition}[theorem] {Proposition}  
\newtheorem{fact}[theorem]{Fact}

\newtheorem*{claim*}{Claim}
\newtheorem*{fact*}{Fact} 
\newtheorem{remark}[theorem] {Remark} 
\newtheorem{definition}[theorem] {Definition} 
\newtheorem*{remark*} {Remark}

\newcommand{\By}[2]{\overset{\mbox{\tiny{#1}}}{#2}}
\newcommand{\ByRef}[2]{   \By{\eqref{#1}}{#2} }

\newcommand{\leBy}[1]{    \By{#1}{\le} }

\newcommand{\eqByRef}[1]{ \ByRef{#1}{=} }

\newcommand{\leByRef}[1]{ \ByRef{#1}{\le} }

\renewcommand{\leq}{\leqslant}
\renewcommand{\le}{\leqslant}

\renewcommand{\ge}{\geqslant}

\renewcommand{\epsilon}{\varepsilon}

\title{Connectivity of inhomogeneous random graphs II}
\author{Jan Hladk\'y}
\address{Institute of Computer Science of the Czech Academy of Sciences, Pod Vod\'{a}renskou v\v{e}\v{z}\'{\i} 2, 182~07 Prague, Czechia. With institutional support RVO:67985807.}
\email{hladky@cs.cas.cz}
\author{Gopal Viswanathan}
\address{\emph{This work was done while affiliated at:} Institute of Computer Science of the Czech Academy of Sciences, Pod Vod\'{a}renskou v\v{e}\v{z}\'{\i} 2, 182~07 Prague, Czechia. With institutional support RVO:67985807. \emph{and} Charles University, Faculty of Mathematics and Physics, Czechia}
\thanks{Research supported by Czech Science Foundation Project 21-21762X}
\subjclass[2010]{05C80}
\keywords{random graph, graphon, connectedness, isolated vertex, minimum degree}
\date{}

\begin{document}
\begin{abstract}
Each graphon $W:\Omega^2\rightarrow[0,1]$ yields an inhomogeneous random graph model $\RG(n,W)$. We show that $\RG(n,W)$ is asymptotically almost surely connected if and only if \emph{(i)} $W$ is a connected graphon \emph{and} \emph{(ii)} the measure of elements of $\Omega$ of $W$-degree less than $\alpha$ is $o(\alpha)$ as $\alpha\rightarrow 0$. These two conditions encapsulate the absence of several linear-sized components, and of isolated vertices, respectively.

We study in bigger detail the limit probability of the property that $\RG(n,W)$ contains an isolated vertex, and, more generally, the limit distribution of the minimum degree of $\RG(n,W)$.
\end{abstract}
\maketitle

\section{Introduction}
This paper studies the connectivity of graphon-based inhomogeneous random graphs. 
We begin by recalling results about connectivity in the Erdős--Rényi random graph. This is a classical topic, and indeed was the subject of the first two (and independent) papers on random graphs by Gilbert~\cite{MR108839} and by Erdős and Rényi~\cite{MR120167} in~1959. Their results imply that for each $\eps>0$ the random graph $\RG(n,(1-\eps)\frac{\ln n}{n})$ is asymptotically almost surely disconnected while $\RG(n,(1+\eps)\frac{\ln n}{n})$ is asymptotically almost surely connected. There is a specific reason for disconnectedness: $\RG(n,(1-\eps)\frac{\ln n}{n})$ asymptotically almost surely contains an isolated vertex. In fact, an even stronger connection between connectedness and isolated vertices is known (see e.g.\ Theorem~4.2 in~\cite{MR3675279}) concerning what is often called the Erdős--Rényi graph process: starting with an $n$-vertex edgeless graph and adding uniformly random edges one-by-one, asymptotically almost surely, the graph becomes connected exactly at the point when the last isolated vertex disappears. Similar local obstructions govern even more complicated connectivity-related objects. For example, famous results of~\cite{HittingBol} and of~\cite{HittingAKSz} say that in the above Erdős--Rényi graph process asymptotically almost surely, the graph becomes Hamiltonian exactly at the point when the last vertex  of degree~1 disappears.

In this paper, we look at graphon-based inhomogeneous random graphs. A \emph{graphon} is a symmetric measurable function $W:\Omega^2\rightarrow [0,1]$, where $(\Omega,\mu)$ is a probability space with an implicit sigma-algebra. Each graphon $W$ yields for each $n\in \NN$ a random graph model $\RG(n,W)$ defined as follows. The vertex set of $G\sim \RG(n,W)$ is the set $[n]$. To generate edges of $G$, first, sample $n$ independent elements $x_1,\ldots,x_n$ in $\Omega$ according to $\mu$. Second, for each pair $ij$ independently include $ij$ as an edge of $G$ with probability $W(x_i,x_j)$. The first stage will be referred to as the \emph{vertex-generating stage} and the second as the \emph{edge-generating stage}.

Graphon-based inhomogeneous random graphs generalize Erdős--Rényi random graphs $\RG(n,p)$ when $p\in[0,1]$ is constant. To this end it is enough to consider the graphon which is constant $p$. In the 1980s, motivated by social networks, the so-called \emph{stochastic block model} was proposed, see~\cite{MR718088}. The stochastic block model (roughly) corresponds to the case when the graphon $W$ in $\RG(n,W)$ is a step-function on an $N\times N$-grid (for some fixed $N\in\NN$) on $\Omega^2$, thus being half-way in generality and complexity between Erdős--Rényi random graphs $\RG(n,p)$ and the full inhomogeneous model. The model $\RG(n,W)$ in its full generality emerged in~2004 as one of the key concepts of the then-surfacing theory of dense graph limits. Lovász and Szegedy~\cite{Lovasz2006} proved that each sequence of graphs contains a subsequence which in an appropriate sense converges to a graphon. This in itself does not justify graphons as the ultimate limits of graphs as it could be that not all graphons are needed for this representation. This is where inhomogeneous random graphs enter the picture. Indeed, one of the key results in~\cite{Lovasz2006} is that for each graphon $W$ a sequence $(\RG(n,W))_n$ converges almost surely to $W$. In other words, each graphon is needed.

The Erdős--Rényi random graph $\RG(n,p)$ when $p$ is constant corresponds to the dense case in which the random graph typically has $\Theta(n^2)$ edges. This is the most tractable case and challenges (many of them still open) concerning parameters such as the independence number, the chromatic number, or features such as tree universality arise in much sparse regimes such as $p=p(n)=\Theta(\ln n/n)$ or $p=p(n)=\Theta(1/n)$. Excluding the trivial case $W\equiv 0$, the random graph $\RG(n,W)$ is dense. However, the intuition that dense is easy does not carry over to $\RG(n,W)$. The reason is obviously in the inhomogeneities of $W$. A particular source of difficulties in behavior of $\RG(n,W)$ are areas of $\Omega^2$ in which $W(x,y)$ approaches~$0$. In particular, these areas cause that parts of $\RG(n,W)$ behave like a sparse random graph, a phenomenon that plays an important role in the problem we study here.

Since then, various properties of inhomogeneous random graphs have been studied. This includes results on the clique number~\cite{MR3683364,MR4011679}, the diameter~\cite{MR3845099}, and the chromatic number~\cite{IK:ChromaticGraphon}.  

\subsection{Our results}
Our results concern the minimum degree and the vertex connectivity of $\RG(n,W)$. The minimum degree of a graph $G$ is denoted by $\delta(G)$. Recall that the \emph{vertex connectivity} of $G$, denoted $\kappa(G)$, is the minimum number of vertices whose deletion from $G$ results in a disconnected graph or in a graph with at most one vertex; in particular, $\kappa(K_n)=n-1$ for the complete graph on $n$ vertices. Throughout the paper we refer to $\kappa(G)$ simply as the \emph{connectivity} of $G$. These two parameters are related by the trivial general inequality $\delta(G)\ge \kappa(G)$.

We start by discussing the minimum degree of $\RG(n,W)$. To this end recall the notion of degrees in graphons. For $x\in \Omega$, define
\begin{equation}\label{eq:degrees}
\deg_W(x):=\int_{\Omega}W(x,y)\differential\mu(y).
\end{equation}
We have that for $\mu$-almost every $x$, $\deg_W(x)$ is well-defined and is between~0 and~1. We open by the following easy fact about the minimum degree of $\RG(n,W)$. Here, we recall that $\essinf_{x\in\Omega}\deg_W(x)$ is the essential infimum of the degree function, which is defined as $$\essinf_{x\in\Omega}\deg_W(x):=\inf\big\{\gamma>0: \mu(\{x\in\Omega:\deg_W(x)<\gamma\})>0\big\}.$$
\begin{fact}\label{fact:linearmindeg}
Suppose that $W$ is a graphon. Write $\beta:=\essinf_{x\in\Omega}\deg_W(x)$. Then asymptotically almost surely, we have $\delta(\RG(n,W))=(\beta+o(1))n$.
\end{fact}
We could not find a reference for this statement; we sketch the easy proof. Suppose that in the vertex-generating stage we condition that $x_i=x$ for some $x\in \Omega$. Then for each $j\in[n]\setminus \{i\}$ the event $ij\in E(\RG(n,W))$ occurs with probability $\deg_W(x)$, and these events are independent for different $j$'s. 
That is, for $G\sim \RG(n,W)$, the degree $\deg_{\RG(n,W)}(i)$ of vertex $i$ in $G$ has binomial distribution with parameters $n-1$ and $\deg_W(x)$.
Chernoff's bound tells us that we have exponential concentration around
\begin{equation}\label{eq:degreeBin}
\deg_{\RG(n,W)}(i)\overset{\mathrm{distr}}{\approx}(n-1)\deg_W(x)\approx n\deg_W(x).
\end{equation}
Since the concentration is exponential, a union bound over $i\in[n]$ gives that asymptotically almost surely $\deg_{\RG(n,W)}(i)=(\deg_W(x_i)+o(1))n$ for all $i\in[n]$. To conclude the proof of Fact~\ref{fact:linearmindeg}, we observe that, almost surely, in the vertex-sampling we have $\deg_W(x_i)\ge \beta$ for all $i\in [n]$, and also asymptotically almost surely we have $\deg_W(x_i)\le \beta+o(1)$ for at least one $i\in [n]$.

We want to refine the bound $\delta(\RG(n,W))=(\beta+o(1))n$ in the case when $\essinf_{x\in\Omega}\deg_W(x)=0$. In particular, we want to understand the sparsest possible regime when $\delta(\RG(n,W))$ is constant, which in particular carries the information whether $\RG(n,W)$ contains at least one isolated vertex. It is instructive to consider the following one-parameter family of graphons $(U_t)_{t\in(0,\infty)}$ where $U_t:(0,1]^2\rightarrow [0,1]$ is defined by $U_t(x,y)=x^ty^t$. Thus, for a graphon $U_t$ and for $x\in(0,1]$ we have $\deg_{U_t}(x)=\Theta(x^t)$. After the vertex-generating stage we typically have $\min_{i\in[n]} x_i =\Theta(n^{-1})$. Combined with the correspondence between degrees in the graphon and in the random graph~\eqref{eq:degreeBin}, we have $\delta(\RG(n,U_t))\approx \Theta(n^{-t})n$. In particular, for $t>1$ we expect to typically have $\delta(\RG(n,U_t))=0$; for $t<1$ we expect to typically have that $\delta(\RG(n,U_t))$ is superconstant; and the case $t=1$ seems to lead to a more subtle behavior. All these heuristics could be turned into fairly easy rigorous calculations. This example suggests that the tail behavior of the measure of vertices of small degree in the graphon is a key quantity to look at.
\begin{definition}[Key notation $S_W(\alpha)$ and $g_W(\alpha$)]
Given a graphon $W$ and a number $\alpha\in[0,1]$, define $S_W(\alpha)\subset \Omega$ by Further, define $g_W(\alpha)\in[0,1]$ by $g_W(\alpha):=\mu(S_W(\alpha))$.
\end{definition}
We can now state our main theorem concerning $\delta(\RG(n,W))$.
\begin{theorem}\label{thm:mindegree}
Suppose that $W$ is a graphon.
\begin{enumerate}[label=(\alph*)]
    \item\label{en:iso} If $\lim_{\alpha\searrow0}\frac{g_W(\alpha)}{\alpha}=+\infty$ then a.a.s. $\RG(n,W)$ contains at least one isolated vertex.
    \item  \label{en:iso2} If $\lim_{\alpha\searrow0}\frac{g_W(\alpha)}{\alpha}=0$ then $\delta(\RG(n,W))\to \infty$ in probability.
    \item\label{en:iso3R1} If $\liminf_{\alpha\searrow0}\frac{g_W(\alpha)}{\alpha}>0$ then we have $\liminf_{n\to\infty}\Pr[\textrm{$\RG(n,W)$ contains isolated vertex}]>0$.
    \item\label{en:iso3R2} If $\limsup_{\alpha\searrow0}\frac{g_W(\alpha)}{\alpha}<+\infty$ then for every $\ell\in \NN$ we have $\liminf_{n\to\infty}\Pr[\delta(\RG(n,W))>\ell]>0$.
\end{enumerate}
\end{theorem}
We prove Theorem~\ref{thm:mindegree} in Section~\ref{sec:proofmindegree}. There is one important regime our theorem does not cover, namely the regime when $\lim_{\alpha\searrow0}\frac{g_W(\alpha)}{\alpha}=\beta\in(0,\infty)$. One might expect that in that regime the distribution of $\delta(\RG(n,W))$ converges to a limit distribution which depends only on the parameter~$\beta$. In Section~\ref{ssec:properlimit} we show that this is not the case by constructing two graphons $U_1$ and $U_2$ with $\lim_{\alpha\searrow0}\frac{g_{U_1}(\alpha)}{\alpha}=\lim_{\alpha\searrow0}\frac{g_{U_2}(\alpha)}{\alpha}=\frac23$ where $\delta(\RG(n,U_1))$ converges to the geometric distribution with success parameter $1-\exp(-\frac23)$ and $\delta(\RG(n,U_2))$ converges to the geometric distribution with success parameter $\frac12$.

We now turn to connectivity of $\RG(n,W)$. One obvious obstacle for connectivity is the presence of isolated vertices; we saw that in the case of Erdős--Rényi random graphs it was essentially the only obstacle. In the inhomogeneous setting, there is another macroscopic obstacle. We say that a graphon $W:\Omega^2\rightarrow [0,1]$ is \emph{disconnected} if there exists a set $X\subset \Omega$ with $\mu(X)\in(0,1)$ such that
\begin{equation}\label{eq:disconnected}
	\int_{X\times(\Omega\setminus X)}W(x,y)\differential\mu^2(x,y)=0.
\end{equation}
We say that $W$ is \emph{connected} if it is not disconnected. The concept was first studied by Janson in~\cite{Janson08Connectedness}. Suppose that $W$ is a disconnected graphon and $X\subset \Omega$ is a witness to that end. By the law of large numbers, we will asymptotically almost surely have that in the vertex-generating stage the sets $I:=\{i\in[n]:x_i\in X\}$ and $J:=\{i\in[n]:x_i\in \Omega\setminus X\}$ are nonempty. In the edge-generating stage, we shall almost surely have by~\eqref{eq:disconnected} that $e_{\RG(n,W)}(I,J)=0$, and thus $\RG(n,W)$ is asymptotically almost surely disconnected.

Our main theorem says that when the graphon is connected and when no isolated vertices appear, the corresponding inhomogeneous random graph is indeed asymptotically almost surely connected. Actually, we get a positive result for higher connectivity as well.
\begin{theorem}\label{thm:connected}
	Suppose that $W$ is a connected graphon with $\lim_{\alpha\searrow0}\frac{g_W(\alpha)}{\alpha}=0$. Then for every $K\in \NN$ a.a.s. $\RG(n,W)$ is a $K$-connected graph.
\end{theorem}
We prove Theorem~\ref{thm:connected} in Section~\ref{sec:proofconnected}.

In Section~\ref{sec:concluding} we offer conjectures regarding a more precise quantification of our results and a closer relation between connectivity and isolated vertices.

\subsection{Comparison with~\cite{DevFra:Connectivity}}\label{ssec:Devroy}
Our results are very closely related to paper of Devroye and Fraiman~\cite{DevFra:Connectivity}. A \emph{kernel} is a symmetric nonnegative function $\Gamma\in L^1(\Omega^2)$. The fact that the values of $\Gamma$ are not restricted to $[0,1]$ is sensible when we add rescaling of the edge-insertion probabilities to the random graph generating procedure. In particular, \cite{DevFra:Connectivity} studies rescaling by $\frac{\ln n}{n}$. That is, $\RK(n,\Gamma)$ is a random graph model on vertex set $[n]$ and with the edge set generated as follows. First, sample $n$ independent elements $x_1,\ldots,x_n$ in $\Omega$ according to $\mu$. Second, for each pair $ij$ independently include $ij$ as an edge of $\RK(n,\Gamma)$ with probability $\min\{1,\frac{\ln n}{n}\cdot\Gamma(x_i,x_j)\}$. The choice of the rescaling factor $\frac{\ln n}{n}$ was clearly motivated by the connectivity threshold for Erdős--Rényi random graphs.

For kernels, the notions of connectivity and degrees defined via~\eqref{eq:disconnected} and~\eqref{eq:degrees} apply unchanged. In particular, the degree function $\deg_\Gamma(\cdot)$ is in $L^1(\Omega)$.

The main result of~\cite{DevFra:Connectivity} reads as follows.
\begin{theorem}\label{thm:DevFr}
	Suppose that $\Gamma$ is a connected kernel. Define a function $h:\Omega\rightarrow [0,\infty)$ by $h(x):=\int_\Omega \Gamma(x,y)^2$. Define $\lambda=\essinf_{x\in\Omega}\deg_\Gamma(x)$.
	\begin{enumerate}[label=(\alph*)]
		\item\label{en:DevFrNegative} If $\lambda<1$ and $h\in L^2(\Omega)$ then a.a.s. $\RK(n,\Gamma)$ contains an isolated vertex. In particular, a.a.s. $\RK(n,\Gamma)$ is disconnected.
		\item\label{en:DevFrPositive} If $\lambda>1$, $h\in L^\infty(\Omega)$, and $\Gamma$ is a continuous function then a.a.s. $\RK(n,\Gamma)$ is connected.
	\end{enumerate}
\end{theorem}
The conditions of $h\in L^2(\Omega)$, of $h\in L^\infty(\Omega)$, and of the continuity of $\Gamma$ are certain regularity conditions. As the authors show, some conditions are necessary, but these conditions (and most strikingly that of the continuity of $\Gamma$) are probably not optimal. Putting these technical conditions aside, the most important parameter here is obviously $\lambda$, and the distinction $\lambda<1$ versus $\lambda>1$ plays the same role as our distinction $\lim_{\alpha\searrow0}\frac{g_W(\alpha)}{\alpha}=+\infty$ versus $\lim_{\alpha\searrow0}\frac{g_W(\alpha)}{\alpha}=0$. We want to emphasize these two distinctions carry different information about the kernel/graphon. Indeed, $\lambda>1$ means that the degree of \emph{every} element in $\Gamma$ is big, whereas $\lim_{\alpha\searrow0}\frac{g_W(\alpha)}{\alpha}=0$ permits elements of arbitrarily small degree, provided the measure $g_W(\alpha)$ of such elements shrinks faster than $\alpha$ itself as $\alpha\to0$: low-degree elements may exist, but must become rarer as the degree threshold decreases.

As an extension of one tool we develop in this paper, we can strengthen Theorem~\ref{thm:DevFr}\ref{en:DevFrPositive} to higher connectivity. 
\begin{theorem}\label{thm:DevFrBetter}
Suppose that $\Gamma$ is a kernel. Define a function $h:\Omega\rightarrow [0,\infty)$ by $h(x):=\int_\Omega \Gamma(x,y)^2$. Define $\lambda:=\essinf_{x\in\Omega}\deg_\Gamma(x)$.
If $\lambda>1$, $h\in L^\infty(\Omega)$, and $\Gamma$ is a continuous function then for every $K\in \NN$, $\RK(n,\Gamma)$ is a.a.s. $K$-connected.
\end{theorem}
The proof is given in Section~\ref{sec:proofconnected}.

\subsection*{Acknowledgments}
We thank Yixiang Wang and three anonymous referees for their feedback.

\subsection{Notation}
The natural logarithm is denoted by $\ln(\cdot)$ while $\log(\cdot)$ is the logarithm of base~2.

In~\eqref{eq:degrees}, we introduced $\deg_W(x)$. We extend this notion to the
degree of $x$ relative to a set $A\subset\Omega$,
$$
\deg_W(x,A):=\int_A W(x,y)\,\differential\mu(y),
$$
so that $\deg_W(x)=\deg_W(x,\Omega)$.

\section{Proof of Theorem~\ref{thm:mindegree}}\label{sec:proofmindegree}
\subsection{Proof of Theorem~\ref{thm:mindegree}\ref{en:iso}}
For the proof we need a well-known bound concerning the binomial distribution. Suppose that $X\sim \Bin(n,p)$ and $r<np$. Then (see e.g.~\cite[page 151]{MR0228020}),
\begin{equation}\label{eq:binomialsmall}
	\Pr[X\le r]\le \frac{(n-r)p}{(np-r)^2}.
\end{equation}

With this tool we can turn to the proof of Theorem~\ref{thm:mindegree}\ref{en:iso}. Let $\eps>0$ be arbitrary. To prove the statement, we will show that for every sufficiently large $n$ we have 
\begin{equation}\label{eq:isoIsIn}
\Pr[\textrm{$\RG(n,W)$ contains isolated vertex}]> 1-2\eps.
\end{equation}
Let $q:=\frac{\eps}{6}$, $c:=\frac{24}{\eps^2}$, and $k:=\frac{cq}{3}=\frac{4}{3\eps}$.
By assumption~\ref{en:iso}, there exists $\alpha_0>0$ such that $g_W(\alpha)>c\alpha$ for every $\alpha\in(0,\alpha_0]$. We claim that~\eqref{eq:isoIsIn} holds when $n>\frac{\eps}{6\alpha_0}$.

Fix a set $S\subset S_W(\frac{q}n)$ of measure exactly $g:=c\cdot \frac{q}n$. We introduce two random counting variables, $N$ and $R$. Let $N$ be the number of elements $x_i$ sampled from $S$ in the vertex-generating stage. Let $R$ be the number of ordered pairs $(x_i,x_j)$ such that $x_i\in S$ and $ij\in E(\RG(n,W))$. Then \eqref{eq:isoIsIn} follows once we prove that
\begin{equation}\label{eq:compNR}
\Pr[N>R]\ge 1-2\eps.
\end{equation}

First, we look at the random variable $N$. This random variable has binomial distribution with parameters $n$ and $g$. In particular, its expectation is $n g= cq= 3k$. We can therefore use~\eqref{eq:binomialsmall}, 
$$
\Pr[N\le k]\le \frac{(n-k)g}{(ng-k)^2}\le \frac{ng}{(\frac23ng)^2}= \frac{9}{4cq}<\eps
.$$

Next, we turn our attention to the random variable $R$. By symmetry and linearity of expectation, $\Exp[R]=\binom{n}{2}r$, where $r$ is the probability that in the vertex-generating stage the element $x_1$ is sampled from $S$ and further (in the edge-generating stage) $12$ forms an edge in $\RG(n,W)$. Using the defining property of the set $S_W(\frac{q}n)$, we have
$$r=\int_{x\in S}\int_{y\in\Omega} W(x,y)\differential\mu(y)\differential\mu(x)\le \int_{x\in S}\frac{q}{n}\differential\mu(x)=\frac{g q}n.$$
We get $\Exp[R]\le q g n=c q^2$. By Markov's inequality, $\Pr[R\ge k]\le \frac{c q^2}{k}=\frac{\eps}2$. 
The bounds on $\Pr[N\le k]$ and $\Pr[R\ge k]$ imply~\eqref{eq:compNR}.

\subsection{Proof of Theorem~\ref{thm:mindegree}\ref{en:iso2}}
To prove the statement, we need to prove that for every $\eps>0$ and every $\ell\in\NN$, when $n$ is sufficiently large, then $\Pr\left[ \delta(\RG(n,W))  < \ell\right] \leq \epsilon$.
Set $c\ge 1$ be large enough that
\begin{align}
	\label{eq:cvelke1}
	4(3ct)^\ell<\exp(ct/4) \quad \text{for each $t\ge\tfrac18$.}
\end{align}
By the assumption~\ref{en:iso2}, there exists $\alpha_0>0$ so that for all $\alpha\in(0,\alpha_0]$ we have $\frac{4cg_W(\alpha)}{\eps}<\alpha$.
Suppose that $n$ is large enough, in particular that
\begin{equation}
	\label{eq:ans}
\frac{c}{\sqrt{n}}<	\alpha_0.
\end{equation}

Let $\deg(i):=\deg_{\RG(n,W)}(i)$ be a random variable that denotes the degree of vertex $i$ in the graph $\RG(n,W)$ (we use this shorter notation from here on).
By the union bound,
\begin{equation}\label{eq:dotheunionbound}
\Pr\left[ \delta(\RG(n,W))  < \ell\right] \leq \sum_{i=1}^{n} \Pr\left[\deg(i)  < \ell\right]= n \cdot \Pr\left[\deg(1)  < \ell\right].
\end{equation}
We shall prove that $\Pr\left[\deg(1)  < \ell\right] \leq \frac{\epsilon}{n}$
and Theorem~\ref{thm:mindegree}\ref{en:iso2} follows.

Recall that $ x_{1}$ is the element sampled from $\Omega$ used to represent the vertex~1.  For  $k \in\NN$, let $B_{k}$ denote the set $S_W(\frac{c2^{k}}{n} )$ and define $B_0:=\emptyset$.

We have
\begin{align}
\Pr\left[\deg(1)  < \ell\right] =& \sum_{k=1}^{\infty} \Pr\left[x_{1} \in B_{k} \setminus B_{k-1} \right] \cdot  \Pr\left[\deg(1)  < \ell \:|\: x_{1} \in   B_{k} \setminus B_{k-1}\right]
\nonumber
\\ 
\nonumber
=&\Pr\left[x_{1} \in B_{1}\right] \cdot \Pr\left[\deg(1)  < \ell \:|\: x_{1} \in   B_{1}\right]\\
&+\sum_{k=2}^{\infty} \Pr\left[x_{1} \in B_{k} \setminus B_{k-1} \right] \cdot  \Pr\left[\deg(1)  < \ell \:|\: x_{1} \in   B_{k} \setminus B_{k-1}\right]
\label{eq:jnG}
\end{align}
We write $T=\Pr\left[x_{1} \in B_{1}\right] \cdot \Pr\left[\deg(1)  < \ell \:|\: x_{1} \in   B_{1}\right]$. We have
\begin{align*}
T \le 
\Pr\left[x_{1} \in B_{1}\right] = g_W(\tfrac{2c}n)\le  \frac{\epsilon}{2n},
\end{align*}
hence  having obtained a satisfactory bound on $T$. We proceed by obtaining a bound on the remaining sum~\eqref{eq:jnG}.

Note that in the generating procedure of $\RG(n,W)$, once the element $x_1$ was placed, the degree $\deg(1)$ has binomial distribution with parameters $n-1$ and $d:=\deg_W(x_1)$. For any $k\ge 2$, in the conditional space $\{x_1\in   B_{k} \setminus B_{k-1}\}$ we have $d\ge \frac{c2^{k-1}}{n}$. Thus,
\begin{align}
\Pr\left[\deg(1)  < \ell \:|\: x_{1} \in   B_{k} \setminus B_{k-1}\right]
&\le
\sum_{i=0}^{\ell-1} \binom{n-1}{i}\cdot \left(\frac{c2^{k-1}}{n}\right)^i \cdot
\left(1-\frac{c2^{k-1}}{n}\right)^{n-1-i}
\nonumber
\\
&\le \sum_{i=0}^{\ell-1} \left(\frac{ne}{i}\right)^i \cdot \left(\frac{c2^{k-1}}{n}\right)^i \cdot 2
\left(1-\frac{c2^{k-1}}{n}\right)^{n-i}
\nonumber
\\
\JUSTIFY{$i \le n/2$}&\le \sum_{i=0}^{\ell-1} \left(\frac{ec2^{k-1}}{i}\right)^i \cdot 2
\exp(-c2^{k-2})
\le \left(ec2^{k-1}\right)^\ell \cdot 2
\exp(-c2^{k-2})
\nonumber
\\
\JUSTIFY{by~\eqref{eq:cvelke1}}&\le \exp(-c2^{k-3}).
\label{eq:Misa}
\end{align}
We plug this bound into~\eqref{eq:jnG}. For terms $k\ge \frac{1}{2}\log n$ we use trivial bounds $g_W(\tfrac{c2^{k}}{n})\le 1$ and for $k\ge \log n$ also $(g_W(\tfrac{c2^{k}}{n})-g_W(\tfrac{c2^{k-1}}{n}))=0$. For terms $k< \frac{1}{2}\log n$ we use that $g_W(\tfrac{c2^{k}}{n})-g_W(\tfrac{c2^{k-1}}{n})\le g_W(\tfrac{c2^{k}}{n})\le \frac{2^{k}}{n}$ by~\eqref{eq:ans}.
\begin{align*}
\Pr\left[\deg(1)  < \ell\right]-T&\le \sum_{k=2}^{\lfloor\frac{1}{2}\log n\rfloor}\frac{2^{k}}{n}\cdot  \Pr\left[\deg(1)  < \ell \:|\: x_{1} \in   B_{k} \setminus B_{k-1}\right]\\
&~~+\sum_{k=\lceil\frac{1}{2}\log n\rceil}^{\lceil\log n\rceil} \Pr\left[\deg(1)  < \ell \:|\: x_{1} \in   B_{k} \setminus B_{k-1}\right]\\
\JUSTIFY{by~\eqref{eq:Misa}}&\le \sum_{k=2}^{\lfloor\frac{1}{2}\log n\rfloor}\frac{2^{k}}{n}\cdot \exp(-c2^{k-3})+\sum_{k=\lceil\frac{1}{2}\log n\rceil}^{\lceil\log n\rceil}\exp(-c2^{k-3}).
\end{align*}
Now, the first sum is at most $\frac{\eps}{4n}$ by~\eqref{eq:cvelke1}. As for the second sum, observe that its first summand is at most $\exp(-c\sqrt{n}/8)\ll \frac{1}{n^2}$, and this summand is the largest one, and further that there are less than $\eps n/4$ summands.

\subsection{Proof of Theorem~\ref{thm:mindegree}\ref{en:iso3R1}}This proof follows by modifying the proof of Theorem~\ref{thm:mindegree}\ref{en:iso}.
Let $\zeta:=\liminf_{\alpha\searrow0}\frac{g_W(\alpha)}{\alpha}$ and let $n_0$ be such that for every $n\ge n_0$ we have $g_W(\frac{\zeta}{n})>\frac{\zeta^2}{2n}$.

Fix a set $S\subset S_W(\frac{\zeta}{n})$ of measure exactly $g:=\frac{\zeta^2}{2n}$. Let $N$ and $R$ be defined as in~\ref{en:iso}. Many of the calculations below are brief, only capturing differences to~\ref{en:iso}.

The random variable $N$ has binomial distribution with parameters $n$ and $g$. Hence,
$$
\Pr[N\ge 1]=1-\left(1-g\right)^n\ge 1-\exp(-\zeta^2/2)\ge \zeta^2/2
.$$

For the random variable $R$ we have $\Exp[R]=\binom{n}{2}r$, where
$$r=\int_{x\in S}\int_{y\in\Omega} W(x,y)\le \int_{x\in S}\frac{\zeta}{n}=\frac{\zeta^3}{2n^2}.$$
We get $\Exp[R]\le \zeta^3/4$. By Markov's inequality, $\Pr[R\ge 1]\le \zeta^3/4$. We have an isolated vertex whenever the event $\{N\ge 1\}\setminus \{R\ge 1\}$ occurs. Hence, $$\Pr[\textrm{$\RG(n,W)$ contains isolated vertex}]\ge \zeta^2/2-\zeta^3/4>0,$$
as was needed.

\subsection{Proof of Theorem~\ref{thm:mindegree}\ref{en:iso3R2}}
Suppose that $W$ and $\ell$ are given. Let $\psi:=\limsup_{\alpha\searrow0}\frac{g_W(\alpha)}{\alpha}$. Let $\alpha_0>0$ be such that for all $\alpha\in(0,\alpha_0)$ we have $\frac{g_W(\alpha)}{\alpha}<2\psi$.
Set $c\ge 1$ be large enough that
\begin{align}
	\label{eq:cvelke13}
	4\cdot \max(1,\psi)\cdot (3ct)^\ell<\exp(ct/32) \quad \text{for each $t\ge\tfrac18$.}
\end{align}
Suppose that $n$ is large enough, in particular that
\begin{equation}
	\label{eq:ansX}
\frac{c}{\sqrt{n}}<	\alpha_0 \quad\mbox{and}\quad 24c\psi<n.
\end{equation}
For  $k \in\NN$, let $B_{k}:=S_W(\frac{c2^{k}}{n})$. In Lemma~\ref{lem:Z} below we give a technical construction of a certain set $Z$. Let us explain why that set is needed. The task is to prove that with positive probability\footnote{In this outline, \emph{positive probability} is meant uniformly in $n$.} the minimum degree of $\RG(n,W)$ is bigger than a constant $\ell$. If we look at the graphon degrees only through the lenses of the partition $\{B_{k}\}_{k\in\NN}$ of $\Omega$, then surely we cannot deduce any favorable outcome if for at least one element we have $x_i\in B_1$ as that might (in the worst case) mean that $\deg_W(x_i)=0$, being an apparent obstacle for high minimum degree. As $B_1$ is small, with positive probability it occurs that all the vertices are sampled outside of $B_1$. So, we might want to prove that with positive probability $\delta(\RG(n,W))>\ell$ in the conditional space in which all the vertices are sampled outside of $B_1$. The next step of the argument would then be that if $x_i\in B_k\setminus B_{k-1}$ for some $k\ge 2$ then $\deg_{\RG(n,W)}(i)$ has binomial distribution with parameters $n-1$ and $p:=\deg_W(x_i)> \frac{c2^{k-1}}{n}$. This however is not the case because of our conditioning (which applies to all the sampled vertices). With that correction the probability parameter of the binomial distribution is $p:=\frac{\deg_W(x_i,\Omega\setminus B_1)}{1-\mu(B_1)}$, for which again no nontrivial lower-bound exists. That is, there may be elements of $\Omega$ of fairly high degree, but with the unfortunate property that all or almost all of their degree goes to $B_1$. This issue is resolved with the set $Z$ below.
\begin{lemma}\label{lem:Z}
	There exists a set $Z\subset\Omega$ such that
	\begin{enumerate}[label=(\alph*)]
		\item\label{en:B1Z} $B_1\subset Z$,
		\item\label{en:muZ} $\mu(Z)<\frac{12c\psi}{n}$,  and 
		\item\label{en:degdeg} for every $x\in\Omega\setminus Z$ we have $\deg_W(x,\Omega\setminus Z)\ge \frac14\deg_W(x)$.
	\end{enumerate}
\end{lemma}
Before proving Lemma~\ref{lem:Z}, we need to introduce an auxiliary theory. Suppose that $H$ is a finite graph for which $V(H)\subset \NN_0$ and $V(H)\neq\emptyset$ equipped with a weight function $w:E(H)\rightarrow [0,\infty)$. We say that the pair $\mathcal{S}=(H,w)$ is a \emph{sprout} if it satisfies the following `flow reduction' condition for every $i\in V(H)\setminus \{\min(V(H))\}$,
\begin{equation}\label{eq:FR}
\sum_{j\in  V(H)\cap [0,i-1]: \{i,j\}\in E(H)} w(\{i,j\})\ge 3\sum_{j\in  V(H)\cap [i+1,\infty): \{i,j\}\in E(H)} w(\{i,j\}).
\end{equation}
The number $o(\mathcal{S}):=\min(V(H))$ is the \emph{origin of $\mathcal{S}$}. The number $f(\mathcal{S}):=3\sum_{e\in E(H):e\ni o(\mathcal{S})}w(e)$ is the \emph{flow into $\mathcal{S}$}. The number $t(\mathcal{S}):=\sum_{e\in E(H)}w(e)$ is the \emph{total weight of $\mathcal{S}$}. The number $d(\mathcal{S}):=\max (V(H))-o(\mathcal{S})$ is the \emph{depth of $\mathcal{S}$}.

\begin{figure}
	\centering
	\includegraphics[scale=0.8]{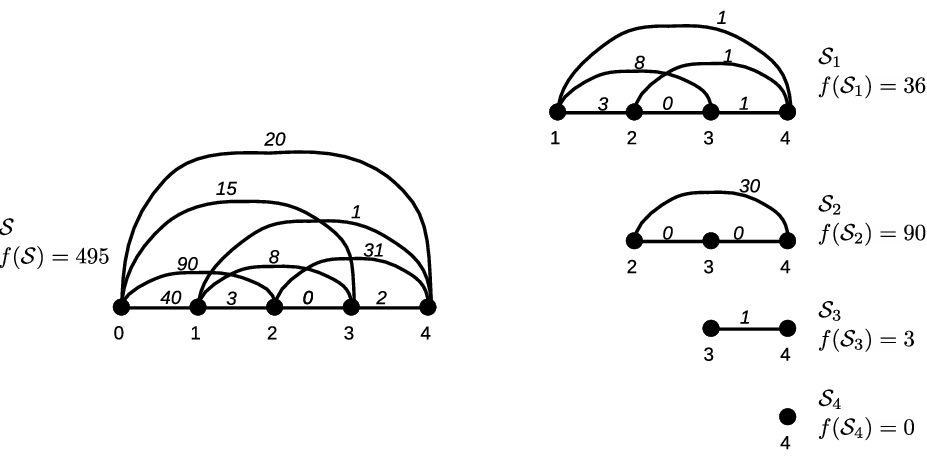}
	\caption{An example of a sprout $\mathcal{S}$ and a sprout family $\{\mathcal{S}_i\}_{i\in[4]}$ covering it.}
	\label{fig:sprouts} 
\end{figure}

Suppose that $\mathcal{S}=(H,w)$ is a sprout, and that $\{\mathcal{S}_\ell=(H_\ell,w_\ell)\}_{\ell\in V(H): \{o(\mathcal{S}),\ell\}\in E(H)}$ is a family of sprouts. We that $\{\mathcal{S}_\ell\}_\ell$ \emph{cover} $\mathcal{S}$ if for every $\ell\in V(H)$ with $\{o(\mathcal{S}),\ell\}\in E(H)$ we have that $H_\ell\subset H$, $o(\mathcal{S}_\ell)=\ell$ and
\begin{equation}\label{eq:flowcover}
f(\mathcal{S}_\ell)\le w(\{o(\mathcal{S}),\ell\}),
\end{equation} 
and further for every $e\in E(H)$ that is not incident with $o(\mathcal{S})$ we have $w(e)=\sum_{\ell: e\in E(H_\ell)}w_\ell(e)$. An example is given in Figure~\ref{fig:sprouts}.

We shall now state two lemmas about sprouts. Lemma~\ref{lem:sproutdecompose} is a preparation for Lemma~\ref{lem:sproutweight} which plays a key role in the proof of Lemma~\ref{lem:Z}.
\begin{lemma}\label{lem:sproutdecompose}
Let $\mathcal{S}=(H,w)$ be an arbitrary sprout. Then there exists a family $\{\mathcal{S}_\ell\}_{\ell\in V(H): \{o(\mathcal{S}),\ell\}\in E(H)}$ that is a cover of $\mathcal{S}$.
\end{lemma}
\begin{proof}
We prove the claim by induction on the number of edges of the sprout carrying nonzero weight. The claim is trivial when this number is~0.

Suppose now that $\mathcal{S}=(H,w)$ is a sprout with at least~1 edge carrying positive weight. In particular, there exists $h\in V(H)$ such that $\{o(\mathcal{S}),h\}$ is an edge of $H$ with $w(\{o(\mathcal{S}),h\})>0$. Let us take $h$ minimal possible. We shall construct two sprouts $\mathcal{S}^-=(H^-,w^-)$ and $\mathcal{S}_h=(H_h,w_h)$ such that $H^-=H-\{h\}$, $o(H_h)=h$, $H_h\subset H$ and for each $e\in E(H^-)\setminus E(H_h)$ and for each $e'\in E(H^-)\cap E(H_h)$ we have 
\begin{equation}\label{eq:edecomp}
w(e)=w^-(e)\quad \mbox{and}\quad w(e')=w^-(e')+w_h(e').
\end{equation}
Further, we shall have that each edge $\tilde{e}$ of $H$ incident to $h$ is present also in $H_h$ and that $w(\tilde{e})=w_h(\tilde{e})$. This is sufficient to prove the inductive step as the sprout $\mathcal{S}^-$ has one edge less and thus can be decomposed further.

Let us proceed with the construction of $\mathcal{S}^-$ and $\mathcal{S}_h$. Actually, we shall only construct $\mathcal{S}_h$ since $\mathcal{S}^-$ is defined thanks to~\eqref{eq:edecomp} as the complement of $\mathcal{S}_h$ in $\mathcal{S}$. We take $H_h:=H[V(H)\cap [h,\infty)]$. Initially, we declare the vertex $h$ \emph{active} and define its \emph{potential} as $\frac{1}{3}w(\{o(\mathcal{S}),h\})$. Now, suppose that at a certain stage we have an active vertex $i\in V(H_h)$ with potential $p_i$. Let $V_i$ denote all $j\in V(H_h)\cap [i+1,\infty)$ with $\{i,j\}\in E(H_h)$. For all $j\in V_i$ we shall define in a single batch all the weights $w_h(\{i,j\})$ arbitrarily subject to
\begin{enumerate}
	\item[(w1)] $w_h(\{i,j\})\le w(\{i,j\})$, and
	\item[(w2)] $\sum_{j\in V_i}w_h(\{i,j\})=\min\left\{p_i,\sum_{j\in V_i}w(\{i,j\})\right\}$.
\end{enumerate} 
Obviously, such an assignment $\{w_h(\{i,j\})\}_{j\in V_i}$ is possible. Having this done, we finish processing this active vertex $i$, and declare vertex $i':=\min (V(H_h)\cap [i+1,\infty))$ as the new active vertex (the procedure finishes if no such vertex exists). We define the potential $p_{i'}$ of this vertex as $p_{i'}:=\frac13\sum_{j\in V(H_h),j<i',\{j,i'\}\in E(H_h)}w_h(\{j,i'\})$, and repeat this procedure.

It remains to argue that $\mathcal{S}^-$ and $\mathcal{S}_h$ are sprouts, that is that they satisfy the flow reduction condition~\eqref{eq:FR}. For $\mathcal{S}_h$ this follows directly from the way potentials are defined and used in~(w2). Let us now turn to $\mathcal{S}^-$. First, for each $i\in V(H^-)\cap [0,h-1]$, the weights of all the edges incident to $i$ are the same as in $\mathcal{S}$, so the flow reduction condition is inherited from $\mathcal{S}$. We skip the case $i=h$ for now and return to it later. For each $i\in V(H^-)\cap [h+1,\infty)$ we distinguish two cases, depending on whether (*): $p_i<\sum_{j\in V_i}w(\{i,j\})$ or (**): $p_i\ge\sum_{j\in V_i}w(\{i,j\})$. In the former case (*), we have
\begin{align*}
\sum_{j\in  [0,i-1]: \{i,j\}\in E(H)} w^-(\{i,j\})&\eqByRef{eq:edecomp}
\sum_{j<i: \{i,j\}\in E(H)\setminus E(H_h)} w(\{i,j\})
+
\sum_{j<i: \{i,j\}\in E(H_h)} \left(w(\{i,j\})-w_h(\{i,j\})\right)
\\
\JUSTIFY{def of $p_i$}&=
\sum_{j<i: \{i,j\}\in E(H)} w(\{i,j\})-3p_i
\\
\JUSTIFY{\eqref{eq:FR} for $\mathcal{S}$}&\ge
3\sum_{j>i: \{i,j\}\in E(H)} w(\{i,j\})-3p_i
\\
\JUSTIFY{(w2)}&\ge
3\sum_{j>i: \{i,j\}\in E(H)} \left(w(\{i,j\})-w_h(\{i,j\})\right)
=
3\sum_{j>i: \{i,j\}\in E(H)} w^-(\{i,j\})
,
\end{align*}
and so~\eqref{eq:FR} is satisfied.
In the latter case~(**), we have $\sum_{j\in V_i}w_h(\{i,j\})=\sum_{j\in V_i}w(\{i,j\})$. That is,~\eqref{eq:edecomp} tells us that $\sum_{j\in V_i}w^-(\{i,j\})=0$, and hence~\eqref{eq:FR} holds trivially.

We return to the case $h=i$, which is very similar. Again, it is enough to assume~(*).
\begin{align*}
	\sum_{j\in  [0,h-1]: \{h,j\}\in E(H)} w^-(\{h,j\})&=
	\sum_{j:j<h,\{h,j\}\in E(H)} w(\{h,j\})-w(\{o(\mathcal{S}),h\})
	\\
	\JUSTIFY{def of $p_h$}&=
	\sum_{j<h: \{h,j\}\in E(H)} w(\{h,j\})	- 3p_h\\
	\JUSTIFY{\eqref{eq:FR} for $\mathcal{S}$}&\ge
	3\sum_{j>h: \{h,j\}\in E(H)} w(\{h,j\})-3p_h
	\\
	\JUSTIFY{(w2)}&\ge
3\sum_{j>h: \{h,j\}\in E(H)} \left(w(\{h,j\})-w_h(\{h,j\})\right)\\
&	=
	3\sum_{j>h: \{h,j\}\in E(H)} w^-(\{i,j\})
	,
	\end{align*}
and so~\eqref{eq:FR} is satisfied.
\end{proof}

\begin{lemma}\label{lem:sproutweight}
Let $\mathcal{S}$ be an arbitrary sprout. Then we have $t(\mathcal{S})\le \frac{1}{2}f(\mathcal{S})$.
\end{lemma}
\begin{proof}
We prove the lemma by induction on the depth of the sprout. In the base case $d(\mathcal{S})=0$ there are no edges and the claim is trivial. So, assume that $d(\mathcal{S})>0$. We appeal to Lemma~\ref{lem:sproutdecompose} and obtain a cover $\{\mathcal{S}_\ell\}_{\ell\in V(H): \{o(\mathcal{S}),\ell\}\in E(H)}$ of $\mathcal{S}$. Each sprout $\mathcal{S}_\ell$ has a smaller depth than $\mathcal{S}$ and thus the induction hypothesis applies to it,
\begin{equation}\label{eq:indhyp}
t(\mathcal{S}_\ell)\le \frac{1}{2}f(\mathcal{S}_\ell)\leByRef{eq:flowcover}\frac{1}{2}w(\{o(\mathcal{S}),\ell\}).
\end{equation}

We have
\begin{equation*}
t(\mathcal{S})=\sum_{e\in E(H), e\ni o(\mathcal{S})}w(e)+\sum_{e\in E(H), e\not\ni o(\mathcal{S})}w(e).
\end{equation*}
We use the definition of $f(\mathcal{S})$ on the first sum and the definition of the covering property on the second sum. Thus,
\begin{align*}
	t(\mathcal{S})&=
	\frac13 f(\mathcal{S})+\sum_{\ell\in V(H): \{o(\mathcal{S}),\ell\}\in E(H)}t(\mathcal{S}_\ell)
\leByRef{eq:indhyp} 
	\frac13 f(\mathcal{S})+\frac12\sum_{\ell\in V(H): \{o(\mathcal{S}),\ell\}\in E(H)}w(\{o(\mathcal{S}),\ell\})\\
\JUSTIFY{def of $f(\mathcal{S})$}&=
\frac13 f(\mathcal{S})+\frac16 f(\mathcal{S})=\frac12f(\mathcal{S})	
,
\end{align*}
as was needed.
\end{proof}

\begin{proof}[Proof of Lemma~\ref{lem:Z}]
Set $Z_0:=B_1$. For $k=1,2,\ldots$ inductively, define $Z_k:=\{x\in\Omega\setminus (\cup_{j=0}^{k-1}Z_j):\deg_W(x,\Omega\setminus (\cup_{j=0}^{k-1}Z_j))<\frac14\deg_W(x)\}$. Take $Z:=\cup_{j=0}^\infty Z_j$.
Obviously,~\ref{en:B1Z} and~\ref{en:degdeg} are satisfied, so it only remains to prove~\ref{en:muZ}. For distinct $i,j\in \NN_0$, define $w(\{i,j\}):=\int_{Z_i\times Z_j} W$; note this is symmetric in $i,j$ since $W$ is symmetric.
We claim that 
\begin{equation}\label{eq:cominghome}
\sum_{i=0}^\infty\sum_{j=i+1}^\infty w(\{i,j\})\le \frac{12\psi c^2}{n^2}.
\end{equation}
It suffices to prove that for each $K\in\NN$, $\sum_{i=0}^K\sum_{j=i+1}^K w(\{i,j\})\le \frac{12\psi c^2}{n^2}$. To this end, we observe that the complete graph on vertex set $\{0,1,\ldots,K\}$ equipped with the weight function $w$ is a sprout. Further, the flow into that sprout is equal to 
$$
3\sum_{i=1}^K w(\{0,i\})=3 \int_{Z_0\times \cup_{i=1}^K Z_i} W\le 3\int_{Z_0}\deg_W(x)\le 3 g_W(\tfrac{2c}{n})\cdot\frac{2c}{n}\le \frac{24\psi c^2}{n^2}.
$$
Hence, Lemma~\ref{lem:sproutweight} tells us that the total weight of that sprout is at most $\frac{24\psi c^2}{n^2}$ as was needed. This establishes~\eqref{eq:cominghome}.

Take an arbitrary $x\in \bigcup_{k=1}^\infty Z_k$, say $x\in Z_k$ for some $k>0$. The definition of $Z_k$ tells us that $x$ sends at least three quarters of its degree into $\cup_{j=0}^{k-1}Z_j$. Since $x\notin B_1$, the degree of $x$ is more than~$\frac{2c}{n}$. We arrive at
$$
\mu\left(\bigcup_{k=1}^\infty Z_k\right)\le \frac{\sum_{i=0}^\infty\sum_{j=i+1}^\infty w(\{i,j\})}{\frac34\cdot \frac{2c}{n}}\leByRef{eq:cominghome}
\frac{\frac{12\psi c^2}{n^2}}{\tfrac34\cdot \frac{2c}{n}}=\frac{8\psi c}n.$$
Hence $\mu(Z)=\mu(Z_0)+\mu\left(\bigcup_{k=1}^\infty Z_k\right)\le g_W(\frac{2c}{n})+\frac{8\psi c}n\le \frac{12\psi c}n$, as was claimed in~\ref{en:muZ}.
\end{proof}

We now return to the main part of the proof. Let $Z$ be from Lemma~\ref{lem:Z}. Let $\cE$ be the event that all the points used to construct $\RG(n,W)$ are sampled from $\Omega \setminus Z$. We have
\begin{align}
\nonumber
\Pr\left[ \delta(\RG(n,W))  < \ell\right]&=\Pr[\cE] \Pr\left[ \delta(\RG(n,W))  < \ell\:|\:\cE\right]+ (1-\Pr[\cE])  \Pr\left[ \delta(\RG(n,W))  < \ell\:|\:\cE^{C}\right]\\
\nonumber
&\le \Pr[\cE] \Pr\left[ \delta(G(n,W))  < \ell\:|\:\cE\right]+ (1-\Pr[\cE])\\
\label{eq:pmi}
&=1-\Pr[\cE](1-\Pr\left[ \delta(G(n,W))  < \ell\:|\:\cE\right])
.
\end{align}
We shall refer to the generating procedure for $\RG(n,W)$. Let $\deg(i)$ be a random variable that denotes the  degree of vertex $i$ in the graph $ \RG(n,W)$.

The number of points sampled from $Z$ in the vertex-generating stage is binomial with parameters $n$ and $\mu(Z)\le \frac{12c\psi}{n}$. Therefore, we have
\begin{align}\label{eq:su1}
	\begin{split}
	\Pr[\cE]&\ge (1-\tfrac{12c\psi}{n})^n\\
	\JUSTIFY{$1-x\ge\exp(-2x)$ for $x\in[0,\frac12]$, and~\eqref{eq:ansX}}&\ge \exp(-24c\psi).
	\end{split}
\end{align}

Next, we bound $\Pr[\deg(1) \  < \ell\:|\:\cE]$ from above.
We have
\begin{align*}
\Pr[\mbox{deg(1)} \  < \ell\:|\:\cE] = \sum_{k=2}^{\log n} \Pr[x_{1} \in B_{k} \setminus B_{k-1}\:|\:\cE ] \cdot  \Pr[\deg(1)  < \ell \:|\:\cE \ \mbox{and} \ x_{1} \in   B_{k} \setminus B_{k-1}]. 
\end{align*}

In the generating procedure for $\RG(n,W)$ given $\cE$, once $x_1$ is placed, $\deg(1)$ has binomial distribution with parameters $n-1$ and $d:=\frac{\int_{y\in \Omega\setminus Z} W(x_1,y)}{1-\mu(Z)}$. For any $k\ge 2$, in the conditional space $\cE\cap\{x_1\in   B_{k} \setminus B_{k-1}\}$ we have with crucial help of property~\ref{en:degdeg} that $d\ge \int_{y\in \Omega\setminus Z} W(x_1,y)\ge\frac{c2^{k-2}}{n}$. Thus,
\begin{align}
\Pr[\deg(1)  < \ell \:|\: \cE\mbox{ and } x_{1} \in   B_{k} \setminus B_{k-1}]
&\le
\sum_{i=0}^{\ell-1} \binom{n-1}{i}\cdot \left(\frac{c2^{k-2}}{n}\right)^i \cdot
\left(1-\frac{c2^{k-2}}{n}\right)^{n-1-i}
\nonumber
\\
&\le \sum_{i=0}^{\ell-1} \left(\frac{ne}{i}\right)^i \cdot \left(\frac{c2^{k-2}}{n}\right)^i \cdot 2
\left(1-\frac{c2^{k-2}}{n}\right)^{n-i}
\nonumber
\\
\JUSTIFY{$i \le n/2$}&\le \sum_{i=0}^{\ell-1} \left(\frac{ec2^{k-2}}{i}\right)^i \cdot 2
\exp(-c2^{k-2})
\nonumber
\\
&\le \left(ec2^{k-2}\right)^\ell \cdot 2
\exp(-c2^{k-3})
\nonumber
\\
\JUSTIFY{by~\eqref{eq:cvelke13}}&\le \exp(-c2^{k-4}).
\label{eq:Misa1}
\end{align}
We now put this together for all $k\ge 2$. For terms $k\ge \frac{1}{2}\log n$ we use trivial bounds $g_W(\tfrac{c2^{k}}{n})\le 1$ and for $k\ge \log n$ also $(g_W(\tfrac{c2^{k}}{n})-g_W(\tfrac{c2^{k-1}}{n}))=0$. For terms $k< \frac{1}{2}\log n$ we use that $g_W(\tfrac{c2^{k}}{n})-g_W(\tfrac{c2^{k-1}}{n})\le g_W(\tfrac{c2^{k}}{n})\le \frac{2\psi\cdot c2^{k}}{n}$ by~\eqref{eq:ansX}.
\begin{align*}
\Pr\left[\deg(1)  < \ell\:|\:\cE\right]&\le \sum_{k=2}^{\lfloor\frac{1}{2}\log n\rfloor}
\Pr\left[x_{1} \in   B_{k} \setminus B_{k-1}\:|\:\cE\right]\cdot  \Pr\left[\deg(1)  < \ell \:|\: \cE\mbox{ and } x_{1} \in   B_{k} \setminus B_{k-1} \right]\\
&~~+\sum_{k=\lceil\frac{1}{2}\log n\rceil}^{\lceil\log n\rceil} \Pr\left[\deg(1)  < \ell \:|\: \cE\mbox{ and } x_{1} \in   B_{k} \setminus B_{k-1} \right]\\
\JUSTIFY{by~\eqref{eq:Misa1}}&\le \sum_{k=2}^{\lfloor\frac{1}{2}\log n\rfloor}\frac{2\psi\cdot c2^{k}}{n}\cdot \exp(-c2^{k-4})+\sum_{k=\lceil\frac{1}{2}\log n\rceil}^{\lceil\log n\rceil}\exp(-c2^{k-4}).
\end{align*}
Let us focus on the first sum. We use~\eqref{eq:cvelke13} to each summand and get $\sum_{k=2}^{\lfloor\frac{1}{2}\log n\rfloor}\frac{2\psi\cdot c2^{k}}{n}\cdot \exp(-c2^{k-4})\le \sum_{k=2}^{\lfloor\frac{1}{2}\log n\rfloor}\frac{1}{2^k n}\le \frac{1}{3n}$. As for the second sum, observe that its first term is at most $\exp(-c\sqrt{n}/18)\ll \frac{1}{n^2}$, and this term is the largest one, and further that there are less than $\log n$ summands. We conclude that $\Pr\left[\deg(1)  < \ell\:|\:\cE\right]\le \frac1{2n}$. By the union bound,
\begin{equation}\label{eq:su2}
\Pr[\delta(\RG(n,W))  < \ell\:|\:\cE] \leq n \cdot \Pr[\deg(1)  < \ell\:|\:\cE]\le \frac12
.
\end{equation}
We substitute bounds~\eqref{eq:su1} and~\eqref{eq:su2} into~\eqref{eq:pmi}.
\begin{equation*}
	\Pr\left[ \delta(\RG(n,W))  < \ell\right]\le 
	1-\exp(-24c\psi)\left(1-\frac12\right)
	.
\end{equation*}
This completes the proof of Theorem~\ref{thm:mindegree}\ref{en:iso3R2}.

\section{Proof of Theorem~\ref{thm:connected}}\label{sec:proofconnected}
We actually prove Theorem~\ref{thm:connected} only in the case $K=1$ in Proposition~\ref{prop:singleconnected}, and combine it with Lemma~\ref{lem:boostConnectivity} which allows to boost the connectivity to an arbitrary constant.
\begin{proposition}\label{prop:singleconnected}
	Suppose that $W$ is a connected graphon with $\lim_{\alpha\searrow0}\frac{g_W(\alpha)}{\alpha}=0$. Then a.a.s. $\RG(n,W)$ is a connected graph.
\end{proposition}
\begin{proof}
We want to prove that for any given $\eps\in(0,1)$, the probability that $\RG(n,W)$ is disconnected is at most $\eps$ for all sufficiently large $n$. Let 
\begin{equation}\label{eq:C}
C:=100\exp(10\eps^{-2})\eps^{-1}.
\end{equation}
 Let $\beta\in(0,C^{-1}]$ be such that for all $\alpha\in(0,\beta)$ we have $g_W(C\alpha)<\frac{\alpha}{C}$. The existence of $\beta$ is guaranteed by our main assumption $\lim_{\alpha\searrow0}\frac{g_W(\alpha)}{\alpha}=0$. Let $\gamma\in(0,\frac12]$ be such that
\begin{equation}\label{eq:gamma}
(2+\ln\tfrac1\gamma)\gamma<\frac{C\beta}{100}.
\end{equation}
The existence of $\gamma$ follows from the fact that $\lim_{\gamma\searrow0}(2+\ln\tfrac1\gamma)\gamma=0$.

We say that an $n$-vertex graph $H$ is \emph{microdisconnected} if there exists a set $X\subset V(H)$ with $|X|\in[1,\gamma n]$ such that $e_H(X,V(H)\setminus X)=0$. We say that $H$ is \emph{macrodisconnected} if there exists a set $X\subset V(H)$ with $|X|\in[\gamma n,\frac12 n]$ such that $e_H(X,V(H)\setminus X)=0$. The proof now has two completely different parts.
\begin{align}
\label{eq:micro}
\Pr[\mbox{$\RG(n,W)$ is microdisconnected}]&\le \tfrac{\eps}{2}.\\
\label{eq:macro}
\Pr[\mbox{$\RG(n,W)$ is macrodisconnected}]&\le \tfrac{\eps}{2}.
\end{align}

We first prove~\eqref{eq:macro}. This case can be handled by soft tools from the theory of graphons. In particular, the following result of Janson is key.
\begin{proposition}[Theorem~1.10(i) in~\cite{Janson08Connectedness}]\label{prop:connectedgraphon}
	Suppose that $W$ is a connected graphon. For every $\gamma>0$ there exists $\rho>0$ with the following property. Suppose that $H$ is a graph whose cut distance from $W$ is at most $\rho$. Then for every set $X\subset V(H)$ with $|X|\in[\gamma v(H),\frac12 v(H)]$ we have $e_H(X,V(H)\setminus X)>\rho v(H)^2>0$.
\end{proposition}
So, Proposition~\ref{prop:connectedgraphon} (with the same $W$ and $\gamma$) tells us that each graph whose cut distance from $W$ is at most $\rho$ is macroconnected.
In particular, the Second Sampling Lemma (\cite[Lemma~10.16]{Lovasz2012}) tells us that when $n$ is large this property holds with probability at least $1-\frac{\eps}{2}$ for $\RG(n,W)$.
This finishes the proof of~\eqref{eq:macro}.

We proceed to prove~\eqref{eq:micro}, which is computationally harder.
For $k\in \NN_0$ define $B_k:=S_W(\frac{C2^k}{n})$. Set $Q:=\lfloor\log (\beta n)\rfloor$. In particular, for each $k\in\NN_0\cap [0,Q]$ we have
\begin{equation}\label{eq:wehaveaboundonB}
\mu(B_k)=g_W(\tfrac{C2^k}{n})\le \frac{2^k}{Cn}.
\end{equation}
Write $\cI:=(\NN_0\cap [0,Q])\cup\{\infty\}$. Let $\iota:\Omega\rightarrow\cI$ be the function which assigns to an input $x\in\Omega$ the smallest $k\in \cI\cap [0,Q]$ such that $x\in B_k$ or, if $x\notin B_Q$ it sets $\iota(x):=\infty$.

Suppose that $A\subset [n]$ and $\bh\in\cI^A$. Then the event $\cE_{A,\bh}$ in the vertex-generating stage is defined as $\bigcap_{i\in A}\{\iota(x_i)=\bh_i\}$. Since the information about the domain is contained in $\bh$, we abbreviate $\cE_{\bh}:=\cE_{A,\bh}$.

For a set $A\subset [n]$ we write $\cD_A$ for the event that there are no edges between $A$ and $[n]\setminus A$.
\begin{align}
\nonumber
\Pr[\mbox{$\RG(n,W)$ is microdisconnected}]&\le\sum_{k=1}^{\gamma n}\sum_{A\in\binom{[n]}{k}}\Pr[\cD_A]=\sum_{k=1}^{\gamma n}\binom{n}{k}\Pr[\cD_{[k]}]\\
\nonumber
&\le \sum_{k=1}^{\gamma n}\left(\frac{en}k\right)^k\Pr[\cD_{[k]}]\\
\label{eq:returnhere}
&=
\sum_{k=1}^{\gamma n}\left(\frac{en}k\right)^k\sum_{\bh\in\cI^{[k]}}\Pr[\cE_\bh]\Pr[\cD_{[k]}|\cE_\bh].
\end{align}
To obtain a good upper bound on this sum, we need to partition $\cI^{[k]}$ as follows. Let $\cJ_{k,0}\subset \cI^{[k]}$ consist only of the constant-0 function. For $M\in[Q]$, let $\cJ_{k,M}\subset \cI^{[k]}$ consist of all functions $\bh$ for which $M=\max_{i\in[k]}\bh_i$. Last, for a nonempty set $S\subset [k]$ and $M\in[Q]\cup\{0\}$, let $\cJ_{k,\infty,S,M}\subset \cI^{[k]}$ consist of all functions $\bh$ for which $S=\{i\in[k]:\bh_i=\infty\}$ and $M=\max_{i\in[k]\setminus S}\bh_i$. If $S=[k]$ then we use the convention $\max_{i\in[k]\setminus S}\bh_i=0$.

For $\bh\in\cJ_{k,M}$, let $f(\bh):=kM-\sum_{i\in[k]}\bh_i$. So, $f(\bh)$ measures the deficiency of $\bh$ in the $\ell^1$-norm to the constant-$M$ vector. Likewise, for $\bh\in\cJ_{k,\infty,S,M}$, let $f(\bh):=(k-|S|)M-\sum_{i\in[k]\setminus S}\bh_i$. The following trivial (algebraic) inequalities which hold for all $M\in[Q]\cup\{0\}$ and nonempty $S\subset[k]$ will come in handy,
\begin{equation}\label{eq:handy}
\sum_{\bh\in \cJ_{k,M}}2^{-f(\bh)}\le 2^k \quad\mbox{and}\quad \sum_{\bh\in \cJ_{k,\infty,S,M}}2^{-f(\bh)}\le 2^{k-|S|}.
\end{equation}

Let us work towards an upper bound on $\Pr[\cE_\bh]$. For $\cE_\bh$ to occur, we have that for each $i\in [k]$ such that $\bh_i<\infty$ it holds $x_i\in B_{\bh_i}$. In particular, \eqref{eq:wehaveaboundonB} applies. We distinguish three cases. If $\bh\in\cJ_{k,0}$ then
\begin{equation}\label{eq:PrE0}
\Pr[\cE_\bh]\le \frac{1}{(Cn)^k}.
\end{equation}
If $\bh\in\cJ_{k,M}$ then
\begin{equation}\label{eq:PrEM}
	\Pr[\cE_\bh]\le \frac{2^{\sum_i \bh_i}}{(Cn)^k}=\frac{2^{kM}}{(Cn)^k}\cdot 2^{-f(\bh)}.
\end{equation}
Last, if $\bh\in\cJ_{k,\infty,S,M}$ then similarly
\begin{equation}\label{eq:PrEinfty}
	\Pr[\cE_\bh]\le \frac{2^{(k-|S|)M}}{(Cn)^{k-|S|}}\cdot 2^{-f(\bh)}.
\end{equation}

Next, we deal with the terms $\Pr[\cD_{[k]}|\cE_\bh]$. Actually, we first consider a finer conditioning $\Pr[\cD_{[k]}|x_1,\ldots,x_k]$ for $x_1,\ldots,x_k\in\Omega$. Let us suppose that $i_*\in[k]$ is such that $\deg_W(x_{i_*})=\max_{i\in[k]}\deg_W(x_i)$.
\begin{align}
\nonumber
\Pr[\cD_{[k]}|x_1,\ldots,x_k]&=\Exp_{X_{k+1},\ldots,X_n}\left[\prod_{i\in[k]}\prod_{j=k+1}^n(1-W(x_i,X_j))\right]\\
\nonumber
&=\prod_{j=k+1}^n\Exp_{X_{k+1},\ldots,X_n}\left[\prod_{i\in[k]}(1-W(x_i,X_j))\right]\\
\nonumber
&=\left(\Exp_X\left[\prod_{i\in[k]}(1-W(x_i,X))\right]\right)^{n-k}\\
\nonumber
&\le \left(\Exp_X\left[1-W(x_{i_*},X)\right]\right)^{n-k}\\
\label{eq:tired}
\JUSTIFY{$k\le \tfrac{n}{2}$}&\le\left(1-\deg_W(x_{i_*})\right)^{\frac{n}{2}}
\end{align}
In particular, if $M=\iota(x_{i_*})\neq 0$ (but allowing the case $M=\infty$) then we have $\deg_W(x_{i_*})\ge\frac{C2^{\min(M,Q)}}{2n}$ 
and the bound $1-z\le \exp(-z)$ substituted in~\eqref{eq:tired} gives $\Pr[\cD_{[k]}|x_1,\ldots,x_k]\le \exp(-\frac{C2^{\min(M,Q)}}4)$.
We get that for each $\bh\in\cJ_{k,M}$ we have 
\begin{equation}\label{eq:paprikyM}
\Pr[\cD_{[k]}|\cE_\bh]\le \exp\left(-\frac{C2^{M}}4\right),
\end{equation}
and for each $\bh\in\cJ_{k,\infty,S,M}$,
\begin{equation}\label{eq:paprikyInfty}
	\Pr[\cD_{[k]}|\cE_\bh]\le \exp\left(-\frac{C2^{Q}}4\right)\le \exp\left(-\frac{C\beta n}8\right),
\end{equation}

We now return to~\eqref{eq:returnhere}, splitting the sum $\bh\in\cI^{[k]}$ into one summand $\bh\in\cJ_{k,0}$, and two groups, $\bh\in\cJ_{k,M}$ ($M=1,\ldots,Q)$ and $\bh\in\cJ_{k,\infty,S,M}$ ($S\subset [k]$ nonempty, $M=0,1,\ldots,Q)$. The corresponding bounds are obtained in~\eqref{eq:Sigma1}, \eqref{eq:Sigma2} and \eqref{eq:Sigma3} below. 

First, we deal with $\bh\in\cJ_{k,0}$ in~\eqref{eq:returnhere}. Recall that $\cJ_{k,0}$ is a singleton.
\begin{equation}\label{eq:Sigma1}
	\sum_{k=1}^{\gamma n}\left(\frac{en}k\right)^k\sum_{\bh\in\cJ_{k,0}}\Pr[\cE_\bh]\Pr[\cD_{[k]}|\cE_\bh]\le
	\sum_{k=1}^{\gamma n}\left(\frac{en}k\right)^k\sum_{\bh\in\cJ_{k,0}}\Pr[\cE_\bh]\leByRef{eq:PrE0}	\sum_{k=1}^{\gamma n}\left(\frac{e}{Ck}\right)^k\le \frac{\eps}{10}.
\end{equation}

Next, we deal with $\bh\in\bigcup_{M=1}^Q\cJ_{k,M}$.
\begin{align}
\nonumber
\sum_{k=1}^{\gamma n}\left(\frac{en}k\right)^k\sum_{M=1}^{Q}&\sum_{\bh\in\cJ_{k,M}}\Pr[\cE_\bh]\Pr[\cD_{[k]}|\cE_\bh]\\
\nonumber
&\leBy{\eqref{eq:PrEM},\eqref{eq:paprikyM}}
\sum_{k=1}^{\gamma n}\left(\frac{e}{Ck}\right)^k\sum_{M=1}^{Q}2^{kM}\sum_{\bh\in\cJ_{k,M}}2^{-f(\bh)}\exp\left(-\frac{C2^{M}}4\right)\\
\label{eq:putt}
&\leByRef{eq:handy}
\sum_{k=1}^{\gamma n}\sum_{M=1}^{Q}\left(\frac{2e}{Ck}\right)^k2^{kM}\exp\left(-\frac{C2^{M}}4\right).
\end{align}
We rewrite the individual terms in the exponential form, $$\left(\frac{2e}{Ck}\right)^k2^{kM}\exp\left(-\frac{C2^{M}}4\right)
=
\exp\left(k\ln(2e/C)-k\ln k+kM\ln 2-\frac{C2^M}4\right).$$
If $M\le \frac1{\ln 2}\ln k$ then $-k\ln k+kM\ln 2\le 0$, and so,
\begin{equation}\label{eq:Idfol}
\left(\frac{2e}{Ck}\right)^k2^{kM}\exp\left(-\frac{C2^{M}}4\right)
\le
\exp\left(k\ln(2e/C)-\frac{C2^M}4\right)
\leByRef{eq:C}
\exp\left((k+M)\ln(\eps/30)\right)
.
\end{equation}
If $M>\frac1{\ln 2}\ln k$ then, writing $M':=M-\frac1{\ln 2}\ln k$, we have
\begin{align*}
-k\ln k+kM\ln 2-\frac{C2^{M}}4
&=\ln 2\cdot k\left(M-\tfrac1{\ln 2}\ln k\right)-\frac{Ck2^{M-\frac1{\ln 2}\ln k}}4\\
&=k\left(M'\cdot \ln 2-\frac{C}42^{M'}\right)\leByRef{eq:C} k(1+\lfloor M'\rfloor)\ln(\eps/100),
\end{align*}
and so,
\begin{equation}\label{eq:Ifyl}
\left(\frac{2e}{Ck}\right)^k2^{kM}\exp\left(-\frac{C2^{M}}4\right)
\le \exp(k(1+\lfloor M'\rfloor)\ln(\eps/100)).
\end{equation}
Note that as $M$ ranges through all numbers bigger than $\frac1{\ln 2}\ln k$, the number $M''=\lfloor M'\rfloor$ ranges through all nonnegative integers.

Plugging~\eqref{eq:Idfol} and~\eqref{eq:Ifyl} in~\eqref{eq:putt}, we get
\begin{align}
	\nonumber
	\sum_{k=1}^{\gamma n}\left(\frac{en}k\right)^k\sum_{M=1}^{Q}&\sum_{\bh\in\cJ_{k,M}}\Pr[\cE_\bh]\Pr[\cD_{[k]}|\cE_\bh]\\
	\nonumber
&\le 	
\sum_{k=1}^{\infty}\sum_{M=1}^{\infty}
\exp\left((k+M)\ln(\eps/30)\right)
+
\sum_{k=1}^{\infty}\sum_{M''=0}^{\infty}
\exp(k(1+M'')\ln(\eps/100))
\\
\label{eq:Sigma2}
&\le\frac{\epsilon}{10}.
\end{align}

The last group to treat in~\eqref{eq:returnhere} is $\bh\in\bigcup_{\emptyset\neq S\subset [k]}\bigcup_{M=0}^Q\cJ_{k,\infty,S,M}$.
\begin{align*}
\sum_{k=1}^{\gamma n}&\left(\frac{en}k\right)^k
\sum_{s=1}^{k}
\sum_{S\in\binom{[k]}{s}}
\sum_{M=0}^{Q}\sum_{\bh\in\cJ_{k,\infty,S,M}}\Pr[\cE_\bh]\Pr[\cD_{[k]}|\cE_\bh]\\
\JUSTIFY{wlog $S=[s]$}	&\le \sum_{k=1}^{\gamma n}\left(\frac{en}k\right)^k
	\sum_{s=1}^{k}
	\left(\frac{ek}{s}\right)^s
	\sum_{M=0}^{Q}\sum_{\bh\in\cJ_{k,\infty,[s],M}}\Pr[\cE_\bh]\Pr[\cD_{[k]}|\cE_\bh]\\
	&\leBy{\eqref{eq:PrEinfty},\eqref{eq:paprikyInfty}}
	\sum_{k=1}^{\gamma n}
\left(\frac{en}{k}\right)^k
\sum_{s=1}^{k}
\left(\frac{ek}{s}\right)^s
\sum_{M=0}^{Q}\frac{2^{(k-s)M}}{(Cn)^{k-s}}\sum_{\bh\in\cJ_{k,\infty,[s],M}}2^{-f(\bh)}\exp\left(-\frac{C\beta n}8\right)\\
	&\leByRef{eq:handy}
\sum_{k=1}^{\gamma n}
\left(\frac{en}{k}\right)^k
\sum_{s=1}^{k}
\left(\frac{ek}{s}\right)^s
\sum_{M=0}^{Q}
\frac{2^{(k-s)M}}{(Cn/2)^{k-s}}\exp\left(-\frac{C\beta n}8\right)\\
&\le
\sum_{k=1}^{\gamma n}
\sum_{s=1}^{k}
\sum_{M=0}^{Q}
\exp\left(
2k-(k-s)(\ln k+\ln\tfrac{C}2)+s(\ln n-\ln s)+(k-s)M\ln 2-\frac{C\beta n}8
\right)\\
&\leByRef{eq:gamma}
(Q+1)
\sum_{k=1}^{\gamma n}
\sum_{s=1}^{k}
\exp\left(2k-(k-s)(\ln k-Q\ln 2+\ln\tfrac{C}2)+s(\ln n-\ln s)-\frac{C\beta n}{10}
\right).
\end{align*}
It is pedestrian to check that for a fixed $k$, the argument in $\exp(\cdot)$ is increasing in $s$ (within the range $s=1,\ldots,k$). Therefore,
\begin{align*}
\sum_{k=1}^{\gamma n}&\left(\frac{en}k\right)^k
\sum_{s=1}^{k}
\sum_{S\in\binom{[k]}{s}}
\sum_{M=0}^{Q}\sum_{\bh\in\cJ_{k,\infty,S,M}}\Pr[\cE_\bh]\Pr[\cD_{[k]}|\cE_\bh]\\
&\le (Q+1)
\sum_{k=1}^{\gamma n}
\exp\left(2k+k(\ln n-\ln k) +\ln k -\frac{C\beta n}{10}
\right)\\
& \le 
\exp\left(\ln (Q+1)+2\gamma n  +\ln (\gamma n)+\gamma n (\ln n-\ln (\gamma n))  -\frac{C\beta n}{10}\right)\\
& =
\exp\left(\ln (Q+1) +\ln (\gamma n)+(2 \gamma  +\gamma \ln (\tfrac{1}{\gamma})) n  -\frac{C\beta n}{10}\right).
\end{align*}
The terms in the argument of $\exp(\cdot)$ are all logarithmic in $n$, except the last two terms $(2\gamma+\gamma \ln (\tfrac{1}{\gamma})) n$ and $\frac{C\beta n}{10}$, which are linear in $n$. Hence, for $n$ large, and using~\eqref{eq:gamma} we have 
\begin{equation}\label{eq:Sigma3}
\sum_{k=1}^{\gamma n}\left(\frac{en}k\right)^k
\sum_{s=1}^{k}
\sum_{S\in\binom{[k]}{s}}
\sum_{M=0}^{Q}\sum_{\bh\in\cJ_{k,\infty,S,M}}\Pr[\cE_\bh]\Pr[\cD_{[k]}|\cE_\bh]
\le 
\exp\left(-\frac{C\beta n}{20}\right)<\frac{\eps}{10}.
\end{equation}
This finishes the proof of~\eqref{eq:micro}, and hence of the theorem.
\end{proof}
Next, we state and prove the lemma which allows us to boost connectivity of $\RG(n,W)$ to an arbitrary constant.
\begin{lemma}\label{lem:boostConnectivity}
	Suppose that $W$ is a graphon with the property that $\RG(n,W)$ is a.a.s. a connected graph. Then for every $K\in\NN$ the graph $\RG(n,W)$ is a.a.s. $K$-connected.
\end{lemma}
\begin{proof}
	We shall prove that for each $n\ge 10K+20$ we have
	$$
	0.4^{K+2}\Pr\left[\text{$\RG(n,W)$ is not $K$-connected}\right]\le 
	\Pr\left[\text{$\RG(\lfloor\tfrac{n}2\rfloor,W)$ is not connected}\right]
	.
	$$
	This proves the statement, since $H\sim \RG(\lfloor\tfrac{n}2\rfloor,W)$ can be generated as follows: In the first step, generate $G\sim\RG(n,W)$. In the second step, pick a uniformly random subset $V$ of size exactly $\lfloor\tfrac{n}2\rfloor$, and identify $H:=G[V]$ in which we relabel each $i$th smallest element of $V$ to $i$ (so that after this relabelling, we have $V(H)=\left[\lfloor\tfrac{n}2\rfloor\right]$).
	
	There are many ways in which one can make precise that picking a random subset of $[n]$ of size exactly  $\lfloor\tfrac{n}2\rfloor$ is similar to including each element of $[n]$ independently with probability $0.5$. We shall use the following particular easy statement to this end.
	\begin{fact*}
		Suppose that $n\ge 10 K+20$. Suppose that $A_1,A_2\subset [n]$ are two disjoint sets of size at most $K$ each. Let $V\subset [n]$ be a uniformly chosen set of size exactly $\lfloor\tfrac{n}2\rfloor$. Then we have
		$$\Pr_V[A_1\subset V \text{ and } A_2\cap V=\emptyset]\ge 0.4^{|A_1|+|A_2|}.$$
	\end{fact*}
	
	Let us now consider the above generating procedure for $\RG(\lfloor\tfrac{n}2\rfloor,W)$. Consider the situation that in the first step, $G\sim\RG(n,W)$ is not $K$-connected. We shall prove that in the second step the produced graph is disconnected with probability at least $0.4^{K+2}$, which will conclude the proof. Using that $G$ is not $K$-connected, there exists a set $A_2\subset [n]$ of size at most $K$ so that $G-A_2$ is disconnected. In particular, this means that there exist two vertices $v_1,v_2\in [n]\setminus A_2$ that lie in different components of $G-A_2$. Set $A_1:=\{v_1,v_2\}$. Applying the above fact, we see that with probability at least $0.4^{|A_2|+2}\ge 0.4^{K+2}$ the graph $H$ from the second step excludes entirely $A_2$, but includes the vertices $v_1$ and $v_2$. In that case, $v_1$ and $v_2$ lie in different components of $H$. So, $H$ is disconnected.
\end{proof}
The same technique allows us to prove Theorem~\ref{thm:DevFrBetter}.
\begin{proof}[Proof of Theorem~\ref{thm:DevFrBetter}]
Suppose that $\Gamma$ and $K$ are given.	
	
We rescale $\Gamma$ by a factor of $\lambda^{-3/4}$, thus getting a kernel $\Gamma':=\lambda^{-3/4}\Gamma$ for which we have $\lambda':=\essinf_{x\in\Omega}\deg_\Gamma'(x)=\sqrt[4]{\lambda}>1$. Thus, the kernel $\Gamma'$ satisfies the assumptions of Theorem~\ref{thm:DevFr}\ref{en:DevFrPositive}, and so $\RK(m,\Gamma')$ is a.a.s. connected, as $m\to\infty$. For $n\in \NN$, write $m(n):=\lfloor \lambda^{-1/4}n\rfloor$.
We shall prove that for each $n$ sufficiently large we have
\begin{equation}\label{eq:Ahaba}
\frac{1}{2}\cdot \lambda^{-2/4}\cdot(1-\lambda^{-1/4})^{K}\Pr\left[\text{$\RK(n,\Gamma)$ is not $K$-connected}\right]\le 
\Pr\left[\text{$\RK( m(n),\Gamma')$ is not connected}\right]
.
\end{equation}
This proves the statement; consider the following two-step procedure. First, generate $G\sim \RK(n,\Gamma)$. Second, generate a uniformly random subset $V$ of the vertex set of $G$ of size $m(n)$, let $H:=G[V]$ and relabel its vertices as in the proof of Lemma~\ref{lem:boostConnectivity}. If $i$ and $j$ are two vertices of $V$ represented by elements $x_i$ and $x_j$ of $\Omega$, then the probability that $ij\in E(H)$ is $\min\{1,\frac{\ln n}{n}\cdot\Gamma(x_i,x_j)\}$ which is at least (for $n$ sufficiently large) $\min\{1,\frac{\ln m(n)}{m(n)}\cdot\Gamma'(x_i,x_j)\}$. Thus $H$ stochastically dominates $\RK( m(n),\Gamma')$. That is, to prove~\eqref{eq:Ahaba}, it is enough to prove that in the setting that after the first step $G$ is not $K$-connected, and hence contains a cut set $A_2$ of size at most $K$ and vertices $v_1,v_2$ in different components of $G-A_2$, then with probability at least $\frac{1}{2}\cdot \lambda^{-2/4}\cdot(1-\lambda^{-1/4})^{K}$ the second step produces a disconnected graph $H$. This follows similarly as in the Fact above; the factor $\lambda^{-2/4}$ stands for the probability of keeping both $v_1$ and $v_2$, the factor $(1-\lambda^{-1/4})^{K}$ stands for the probability of removing each vertex of $A_2$, and the factor $\frac{1}{2}$ provides a cushion to transfer between the uniform model (``random subset of size exactly $m(n)$'') and the binomial model (``include each element independently with probability $\lambda^{-1/4}$'').
\end{proof}

\section{Concluding remarks}\label{sec:concluding}
\subsection{The limit probability when $\lim_{\alpha\searrow0}\frac{g_W(\alpha)}{\alpha}\in (0,\infty)$}\label{ssec:properlimit}
As advertised in the introduction, we construct two graphons $U_1$ and $U_2$ with $\lim_{\alpha\searrow0}\frac{g_{U_1}(\alpha)}{\alpha}=\lim_{\alpha\searrow0}\frac{g_{U_2}(\alpha)}{\alpha}=\frac23$ where $\delta(\RG(n,U_1))$ converges to the geometric distribution with success parameter $1-\exp(-\frac23)$ and $\delta(\RG(n,U_2))$ converges to the geometric distribution with success parameter $\frac12$. We thank Matas Šileikis for suggesting some calculations below. 

\begin{figure}
	\centering
	\includegraphics[scale=0.9]{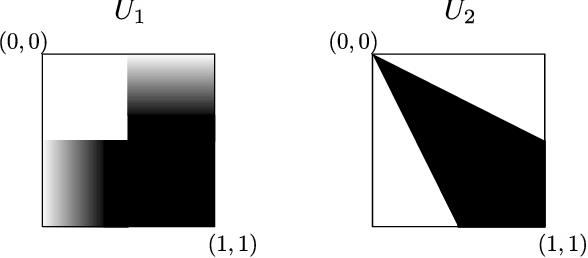}
	\caption{Graphons $U_1$ and $U_2$ from Section~\ref{ssec:properlimit}.}
	\label{fig:twomindegs} 
\end{figure}

Both graphons $U_1,U_2$ are constructed on the domain $[0,1]^2$ and depicted in Figure~\ref{fig:twomindegs}. The graphon $U_1$ is defined to have value~0 on $[0,\frac12)^2$, value~1 on $[\frac12,1]^2$ and value $\min(1,3 y)$ for $y<\frac12$ and $x\ge \frac12$. The graphon~$U_2$ is defined to have value~1 inside a quadrilateral with corners $(0,0)$, $(1,\frac12)$, $(1,1)$ and $(\frac12,1)$ and~0 elsewhere. Indeed, $\lim_{\alpha\searrow0}\frac{g_{U_1}(\alpha)}{\alpha}=\lim_{\alpha\searrow0}\frac{g_{U_2}(\alpha)}{\alpha}=\frac23$.

In the rest of the section we argue what the distribution of the minimum degree in these two models should be. These heuristics could be turned into rigorous arguments in a pedestrian way. In a nutshell, the difference in the distribution of $\delta(\RG(n,U_1))$ and $\delta(\RG(n,U_2))$ reflects that the edges of $\RG(n,U_2)$ are already determined after the vertex-generating stage whereas in $\delta(\RG(n,U_1))$ there is sufficient additional randomness (contributing to the degree of each vertex independently) in the edge-generating stage.

\subsubsection*{Limit distribution of $\delta(\RG(n,U_1))$}
Let us first explain the heuristics for the limit distribution of $\delta(\RG(n,U_1))$. All the above proofs about the minimum degree controlled the number of elements $x_i$ whose degree in the graphon was close to a given quantity. For example in the proof of Theorem~\ref{thm:mindegree}\ref{en:iso2}, this meant controlling how many elements fall into each set $B_k\setminus B_{k-1}$ inside which all the vertices have similar degrees (up to factor~2). A key idea underlying this proof is that only elements $x_i$ of small $\deg_{U_1}(x_i)$ have a non-negligible potential to result in a vertex for which $\deg_{\RG(n,U_1)}(i)=\delta(\RG(n,U_1))$. Here, the relevant scale is $\deg_W(x_i)=\Theta(\frac1n)$. Therefore, for any given $s\in(0,\infty)$, let us look at the number $N_s$ of $i\in [n]$ such that $x_i\in S_{U_1}(\frac{s}n)$, or equivalently, $x_i\le \frac{2s}{3n}$. $N_s$ has binomial distribution with parameters $n$ and $\frac{2s}3$. We conclude that $N_s$ has approximately Poisson distribution $\Poi(\frac23 s)$. In other words, take indices $i_1,i_2,\ldots\in[n]$ such that $x_{i_1},x_{i_2},\ldots$ is the enumeration of all the elements for which we have $\deg_{U_1}(x_{i_j})=\Theta(\frac{1}{n})$, and further such that $\deg_{U_1}(x_{i_1})=\frac{s_1}{n}$, $\deg_{U_1}(x_{i_2})=\frac{s_2}{n}$, \ldots for some (random) constants $0<s_1<s_2<\ldots$. Then the sequence $s_1,s_2,\ldots$ has as similar distribution as a Poisson point process with intensity $\frac23$ on $[0,\infty)$.

If $x_i<\frac12$ and $\deg_{U_1}(x_{i})=\frac{s}{n}$ then the vertex $i$ in $\RG(n,U_1)$ has binomial distribution with parameters $n-1$ and $\frac{s}{n}$. Thus, the distribution of $\deg_{\RG(n,U_1)}(i)$ is approximately $\Poi(s)$. Further, the edge-generating stage is independent for all vertices $i$ with $x_i<\frac12$ (which are the vertices on which the minimum degree of $\RG(n,U_1)$ will be attained almost surely).

Putting this together, we arrive that $\delta(\RG(n,U_1))$ has similar distribution as a random variable $Z$ defined as follows. Let $s_1<s_2<\ldots$ be points generated by the Poisson point process with intensity $\frac23$ on $[0,\infty)$. For each $j\in \NN$ define an independent random variable $Y_j$ with distribution $\Poi(s_j)$. Finally, let $Z:=\min_j Y_j$. In Appendix~\ref{appendix:Minima} we sketch that the distribution of $Z$ is geometric with success parameter $1-\exp(-\frac23)$.

\subsubsection*{Limit distribution of $\delta(\RG(n,U_2))$}
$\delta(\RG(n,U_2))$ is achieved by the degree of vertex $i$ that minimizes the value $x_i$. Further, $i$ is adjacent exactly to vertices $j$ with $x_j\le 2x_i$. That is, $\delta(\RG(n,U_2))$ is determined already in the vertex-generating stage as $\delta(\RG(n,U_2))=|\{j\in [n]: \min_i\{x_i\}<x_j\le 2\min_i\{x_i\}\}|$. With the same reasoning as before, the sequence $x_{i_1},x_{i_2},\ldots$ of the points sampled in the vertex-generating stage, and rescaled by $n$ has approximately the distribution of a Poisson point process with intensity~1.
In Appendix~\ref{appendix:Minima} we show that if $s_1<s_2<\ldots$ are points generated by the Poisson point process with intensity~1 then the number $Z$ of $i>1$ with $s_i<2s_1$ has geometric distribution with success parameter $\frac12$.

\begin{remark}
	Penrose~\cite{MR3476631} studied a model of \emph{soft random geometric graphs}: $n$ points are sampled independently in a Euclidean domain, and each pair of points is joined by an edge with a probability that depends only on the distance between them. This is a natural geometric analogue of our graphon-based model, with Euclidean distance playing the role of the graphon $W$. Despite this similarity, the limiting distribution of the minimum degree found in~\cite{MR3476631} is Poisson, in contrast with the geometric distributions we find for $\delta(\RG(n,U_1))$ and $\delta(\RG(n,U_2))$ in this section. This shows that the geometric limit laws we obtain here are a genuine feature of the graphon setting, rather than a universal phenomenon common to all locally sparse random graph models of this type.
\end{remark}

\subsection{Connectedness versus minimum degree}
The cases $\lim_{\alpha\searrow0} \frac{g_W(\alpha)}{\alpha}=0$ and $\lim_{\alpha\searrow0} \frac{g_W(\alpha)}{\alpha}=\infty$ in Theorems~\ref{thm:mindegree} and Theorems~\ref{thm:connected} show that for a connected graphon $W$ the properties of having an isolated vertex and being disconnected are very tightly related in $\RG(n,W)$. We believe that this relation can be strengthened in two directions at the same time in that
\begin{itemize}
	\item it applies also in the case when the limit $\lim_{\alpha\searrow0} \frac{g_W(\alpha)}{\alpha}$ is in $(0,\infty)$ or it does not exist, and
	\item it puts an asymptotic equivalence between $\delta(\RG(n,W))\ge K$ and $\kappa(\RG(n,W))\ge K$ (as opposed to the simplest case $\delta(\RG(n,W))\ge 1$ and $\kappa(\RG(n,W))\ge 1$).
\end{itemize}
\begin{conjecture}\label{con:onlyoneobstruction}
Suppose that $W$ is a connected graphon. Then for every $K\in \NN$,
$$
\lim_{n\to \infty} \Pr_{G\sim \RG(n,W)}[\delta(G)\ge K\mbox{ and }\kappa(G)<K]=0
.
$$
\end{conjecture}
\begin{figure}
	\centering
	\includegraphics[scale=0.9]{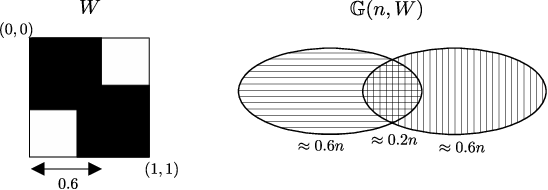}
	\caption{Graphon $W$ and the corresponding random graph $\RG(n,W)$ from Remark~\ref{rem:constobs}.}
	\label{fig:deltakappa} 
\end{figure}

\begin{remark}\label{rem:constobs}
The expressions ``$\delta(G)\ge K\mbox{ and }\kappa(G)<K$'' for each $K$ seem to aggregate to ``$\delta(G)>\kappa(G)$''. However, in Figure~\ref{fig:deltakappa} we show a connected graphon $W$ for which with high probability $\delta(\RG(n,W))\approx 0.6n$ but $\kappa(\RG(n,W))\approx 0.2n$. So, the fact that we control these parameters only in the constant regime is important. (There is a wide window between the constant regime and the linear regime and we do not offer any conjectures for that window.)
\end{remark}

\subsection{A hitting time conjecture}
We mentioned in the introduction that asymptotically almost surely, the hitting time for connectivity equals the hitting time for minimum degree at least~1 in the Erdős--Rényi graph process. Here, we present a counterpart to that conjecture in the inhomogeneous setting, in addition strengthened to higher minimum degree/connectivity as in Conjecture~\ref{con:onlyoneobstruction}.

To this end we present the model, which we call the \emph{$W$-edge-incremental process} of order $n$, and which may be of general interest. This process produces a sequence of random graphs $\mathfrak{X}=(G^t)_{t=0}^{\binom{n}2}$ on vertex set $[n]$, starting with the edgeless graph $G^0$. First, generate elements $x_1,\ldots,x_n\in\Omega$ as in the vertex-generating stage. Next, for each $t=1,\ldots,\binom{n}2$ we construct $G^t$ by adding a single edge to $G^{t-1}$. Call a pair $ij\in\binom{[n]}{2}\setminus E(G^{t-1})$ \emph{visible} if $W(x_i,x_j)>0$, and let $I_t$ be the set of all visible pairs. If $I_t\neq\emptyset$ then pick a random visible pair so that the probability of picking $ab$ is equal to $W(x_a,x_b)/(\sum_{ij\in I_t}W(x_i,x_j))$, and add it to $G^t$. That is, pairs $ab$ with higher values of $W(x_a,x_b)$ tend to be switched to edges earlier than those with smaller values. If $I_t=\emptyset$ then add to $G^t$ a uniformly random nonedge of $G^{t-1}$.\footnote{The construction when $I_t=\emptyset$ is a matter of taste. Alternatively, we could just take $G^t:=G^{t-1}$. Conjecture~\ref{conj:process} is plausible even in this variant.} The construction satisfies the following stochastic dominance sandwiching: $G^{(d-o(1))\binom{n}2}\subset\RG(n,W)\subset G^{(d+o(1))\binom{n}2}$, where $d=\int_{\Omega^2}W$ and the terms $o(1)$ converge to~0 as $n\to\infty$ in probability.
\begin{conjecture}\label{conj:process}
Suppose that $W$ is a connected graphon and $K\in \NN$. For each $n\in\NN$, consider the $W$-edge-incremental process $\mathfrak{X}_n=(G_n^t)_{t=0}^{\binom{n}2}$ of order $n$. Let $f_n$ and $g_n$ be the hitting times for $\mathfrak{X}_n$ for the properties of being $K$-connected and having minimum degree $K$, respectively. Then we asymptotically almost surely have $f_n=g_n$.
\end{conjecture}

\subsection{Hamiltonicity of inhomogeneous random graphs}
Connectivity is a necessary condition for Hamiltonicity, and it is natural to ask whether our
results have a bearing on the (harder) question of when $\RG(n,W)$ is a.a.s.\ Hamiltonian.
After the first version of this paper was written, Garbe, Hladk\'y, and Piga~\cite{GHP:HamiltonicityGnW}
answered this question: they characterized the graphons $W$ for which $\RG(n,W)$ is a.a.s.\
Hamiltonian. Their characterization consists of three conditions on $W$, two of which are
precisely the connectedness of $W$ and the condition $\lim_{\alpha\searrow0}\frac{g_W(\alpha)}{\alpha}=0$
that appear in the hypothesis of Theorem~\ref{thm:connected}. The third condition guarantees
that $\RG(n,W)$ has a.a.s.\ a perfect fractional matching, which is the additional obstruction
one must rule out to go from connectivity to Hamiltonicity.

\appendix
\section{Distributions from Section~\ref{ssec:properlimit}}\label{appendix:Minima}
In this section, we provide calculations for minima of random variables used in Section~\ref{ssec:properlimit}, for which we were unable to find reference in literature. The calculations are simplified, details concerning limit transitions are easy to fill-in. We thank Yixiang Wang for encouraging us to write down these calculations.

Our first proposition is used to justify the distribution of $\delta(\RG(n,U_1))$.
\begin{proposition}
	Suppose that $\lambda>0$ is given. Let $A=\left\{ S_{1},S_{2},\ldots\right\} $
	be a Poisson point process with intensity $\lambda$. Let $\left\{ X_{i}\right\} $
	be a collection of independent random variables where each $X_{i}$
	has Poisson distribution with parameter $S_{i}$. Let $Z=\min_{i}X_{i}$.
	Then $Z$ has geometric distribution with success parameter $1-\exp(-\lambda)$.    
\end{proposition}
\begin{proof}
	Let $m\in\NN_{0}$ be arbitrary. We have to prove that we have
	$\Pr\left[Z>m\right]=\exp(-\lambda (m+1))$.
	
	For $M\in[0,\infty)$, write $Z^{\le M}=\min\left\{ X_{i}:i\in\NN,S_{i}\le M\right\} $.
	Then we have $\Pr\left[Z>m\right]=\lim_{M\to\infty}\Pr\left[Z^{\le M}>m\right]$.
	
	We first express $\Pr\left[Z^{\le M}>m\right]$ for a fixed $M\in (0,\infty)$. Take $L\in\NN$. The little-oh notation refers to $m$ and $M$ as constants, and $L$ going to infinity. For $k\in[L]$, write $I_{M,L,k}$ for the
	interval $I_{M,L,k}=\left(\frac{k-1}{L}\cdot M,\frac{k}{L}\cdot M\right]$.
	We have 
	\[
	\Pr\left[Z^{\le M}>m\right]=\prod_{k=1}^{L}\left(\Pr\left[A\cap I_{M,L,k}=\emptyset\right]+\Pr\left[\left|A\cap I_{M,L,k}\right|=1\right]\cdot p_{M,L,k}\pm\Pr\left[\left|A\cap I_{M,L,k}\right|>1\right]\right),
	\]
	where $p_{M,L,k}=1-\exp\left(-kM/L\right)\sum_{i=0}^{m}\frac{\left(kM/L\right)^{i}}{i!}\pm O(\frac1L)$
	covers the range of probabilities of a $\Poi(x)$-distributed random
	variable to be more than $m$, when $x$ is arbitrary in $I_{M,L,k}$.
	Let us write $b_{M,k,L}=\exp\left(-kM/L\right)\sum_{i=0}^{m}\frac{\left(kM/L\right)^{i}}{i!}$.
	Note also that $\Pr\left[\left|A\cap I_{M,L,k}\right|>1\right]\le(M/L)^{2}$
	for sufficiently large $L$. Hence,
	\begin{align*}
		\Pr\left[Z^{\le M}>m\right] & =\prod_{k=1}^{L}\left(\exp(-\lambda M/L)+\exp(-\lambda M/L)\cdot\frac{\lambda M}{L}\cdot(1-b_{M,L,k})+o\left(\frac{1}{L}\right)\right),\\
		& =\prod_{k=1}^{L}\left((1-\lambda M/L)+(1-\lambda M/L)\cdot\frac{\lambda M}{L}\cdot(1-b_{M,L,k})+o\left(\frac{1}{L}\right)\right),\\
		& =\prod_{k=1}^{L}\left(1-\frac{\lambda M}{L}\cdot b_{M,L,k}+o\left(\frac{1}{L}\right)\right)=\prod_{k=1}^{L}\exp\left(-\frac{\lambda M}{L}\cdot b_{M,L,k}+o\left(\frac{1}{L}\right)\right)\\
		& =\exp\left(-\frac{\lambda M}{L}\sum_{k=1}^{L}b_{M,L,k}\;+\;o(1)\right)\\
		& =\exp\left(-\sum_{i=0}^{m}\frac{1}{i!}\cdot\frac{\lambda M}{L}\sum_{k=1}^{L}\exp\left(-kM/L\right)\left(kM/L\right)^{i}\;+\;o(1)\right).
	\end{align*}
	For any $i=0,\ldots,m$, we look at the summand $\lambda\cdot\frac{M}{L}\sum_{k=1}^{L}\exp\left(-kM/L\right)\left(kM/L\right)^{i}$.
	When $L\to\infty$, this is an approximation of the Riemann integral $\lambda\int_{x=0}^{M}\exp(-x)x^{i}$.
	Recall that $\int_{x=0}^{M}\exp(-x)x^{i}$ is known as the lower incomplete
	gamma function $\gamma(i+1,M)$ and that taking $M\to\infty$, we get
	the gamma function, $\int_{x=0}^{\infty}\exp(-x)x^{i}=i!$. That is
	\begin{align*}
		\Pr\left[Z>m\right] & =\lim_{M\to\infty}\Pr\left[Z^{\le M}>m\right]=\lim_{M\to\infty}\exp\left(-\sum_{i=0}^{m}\frac{1}{i!}\cdot\lambda\cdot\gamma(i+1,M)\right)\\
		& =\exp\left(-\sum_{i=0}^{m}\frac{1}{i!}\cdot\lambda \cdot i!\right)=\exp\left(-\sum_{i=0}^{m}\lambda\right)=\exp\left(-\lambda(m+1)\right),
	\end{align*}
	as was needed.
\end{proof}

The next proposition justifies the distribution of $\delta(\RG(n,U_2))$.
\begin{proposition}
	Suppose that $\lambda>0$ is given. Let $A=\left\{ S_{1},S_{2},\ldots\right\} $
	be a Poisson point process with intensity $\lambda$. Let $B=\min A$ be the left-most point of $A$. Then the number $Z=|\{i\in \NN:S_i\neq B, S_i<2B\}|$ has geometric distribution with success parameter $0.5$ (in particular, it does not depend on $\lambda$).
\end{proposition}
\begin{proof}
	The number $B$ has exponential distribution with intensity parameter $\lambda$, so its probability density function is $x\mapsto\lambda\exp(-\lambda x)$. Given the event $\{B=x\}$, the number of points of $A$ in the interval $(x,2x)$ has Poisson distribution with parameter $\lambda(2x-x)=\lambda x$. We conclude that for each $m\in\NN_0$,
	\begin{align*}
		\Pr[Z=m]&=
		\int_{x=0}^\infty \lambda\exp(-\lambda x)\cdot \frac{\exp(-\lambda x)(\lambda x)^m}{m!}\differential x\\
		&=
		\frac{\lambda}{2^mm!}\int_{x=0}^\infty \exp(-2\lambda x)\cdot (2\lambda x)^m\differential x\\
		\JUSTIFY{substitute $y=2\lambda x$, gamma function}&=0.5^{m+1},
	\end{align*}
	which indeed means that $Z$ has geometric distribution with success parameter~0.5.
\end{proof}

\bibliographystyle{plain}
\bibliography{DSG}

\end{document}